\documentclass[a4paper, 12pt, leqno]{amsart}
\usepackage{mathrsfs}
\usepackage{amsmath}
\usepackage{amscd}
\usepackage{amssymb}
\usepackage{amsthm}
\usepackage{geometry}
\usepackage{tikz}

\title[Based Ring of Two-sided Cells ]
{The Based Rings of Two-sided cells in an Affine Weyl group of type $\tilde B_3$, III}

\author[Y. Qiu]{Yannan Qiu$^{*}$}
\address{$^{*}$
School of mathematics\\
Hangzhou Normal University \\
Zhejiang 311121, China }
 \email{qiuyannan@zju.edu.cn}
\thanks{Partially supported by National Natural Science Foundation of China, No. 12171030.}

\begin{document}
\baselineskip=18pt
\begin{abstract}
We compute the based rings of two-sided cells corresponding to the unipotent class
 in $Sp_6(\mathbb C)$ with Jordan blocks (2211). The results also verify Lusztig's conjecture on the structure of the based rings of the two-sided cells of an affine Weyl group.
\end{abstract}

\maketitle

Let $G$ be $Sp_6(\mathbb C)$ and $W$ be the extended affine Weyl group attached to $G$.  We are concerned with the based rings of two-sided cells of   $W$. In a previous paper [QX2] we computed the  based rings   of the two-sided cells of $W$  with values 2 and 3 of Lusztig's $a$-function.  In this paper we compute the based ring of two-sided cell $c$ of $W$ with value 4 of Lusztig's $a$-function. The result also verifies Lusztig's conjecture on the structure of the based ring of a two-sided cell for $c$. For this two-sided cell, the validity of Lusztig's conjecture on the based rings is already included in main theorem in [BO]. Here we construct the bijection in Lusztig's conjecture explicitly so that the result in this paper can be used for  computing irreducible representations of affine Hecke algebras of type $\tilde B_3$.

The contents of the paper are as follows. Section 1 is devoted to set up and notations. In section 2 we  recall some results on cells of $W$, which are due to J. Du, and establish some conclusions for  the two-sided cell $c$. Sections 3 is devoted to computing the based ring  of the two-sided cell $c$. Under a bijection of Lusztig (see [L3]), this two-side cell corresponds to the unipotent class  in $Sp_6(\mathbb C)$ with Jordan blocks (2211).

\def\Cal{\mathcal}
\def\bold{\mathbf}
\def\ca{\mathcal A}
\def\cdz{\mathcal D_0}
\def\cd{\mathcal D}
\def\cdo{\mathcal D_1}
\def\bold{\mathbf}
\def\l{\lambda}
\def\le{\leq}

\def\ll{\underset {L}{\leq}}
\def\rl{\underset {R}{\leq}}
\def\lr{\rl}
\def\lrl{\underset {LR}{\leq}}
\def\llr{\lrl}
\def\el{\underset {L}{\sim}}
\def\er{\underset {R}{\sim}}
\def\elr{\underset {LR}{\sim}}
\def\ds{\displaystyle\sum}

\section{Set up and notations}

In this section we fix some notations and recall some basic facts. We refer to [KL, L1, L2, L3, QX] for more details.

\medskip

 \noindent{\bf 1.1. The  Hecke algebra of $W$.} \quad Denote by $S$   the set of simple reflections of $W$.
Let $H$ be the Hecke algebra  of $(W,S)$ over $\Cal A=\mathbb
Z[q^{\frac 12},q^{-\frac 12}]$ with
parameter $q$. Then $H$ is a free $\Cal A$-module with standard basis   $\{T_w\}_{w\in W}$. The multiplication in $H$ is defined by $(T_s-q)(T_s+1)=0$ for $s\in S$ and $T_wT_u=T_{wu}$ if $l(wu)=l(w)+l(u)$.
Let
$C_w=q^{-\frac {l(w)}2}\sum_{y\le w}P_{y,w}T_y,\ w\in W$ be the
Kazhdan-Lusztig basis of $H$, where $P_{y,w}$'s are the
Kazhdan-Lusztig polynomials and $\le$ is the Bruhat order of $W$.

We should remark here that our notation $C_w$ stands for $C'_w$ in [KL]. We use the elements $C'_w$ in [KL] since the positivity in multiplications of the elements is simpler in writing and in practical computation. Again for simplicity we write them as $C_w$ in this paper instead of their  original form $C'_w$ in [KL].

The following  formulas for computing $C_w$ (see [KL]) will be used in section 3.

\medskip

\def\vp{\varphi}
\def\st{\stackrel}
\def\sc{\scriptstyle}

(a) Let $w\in W$ and $s\in S$. Then
$$\begin{aligned} C_sC_w=\begin{cases}\displaystyle (q^{\frac12}+q^{-\frac12})C_w,\quad &\text{if\ }sw<w,\\
\displaystyle C_{sw}+\sum_{\st  {z\prec w}{sz<z}}\mu(z,w)C_z,\quad&\text{if\ }sw\ge w,\end{cases}\end{aligned}$$
$$\begin{aligned} C_wC_s=\begin{cases}\displaystyle (q^{\frac12}+q^{-\frac12})C_w,\quad &\text{if\ }ws<w,\\
\displaystyle C_{sw}+\sum_{\st  {z\prec w}{zs<z}}\mu(z,w)C_z,\quad&\text{if\ }ws\ge w,\end{cases}\end{aligned}$$
where $\mu(z,w)$ is the coefficient of the term $q^{\frac12(l(w)-l(z)-1)}$ in $P_{z,w}$ and $z\prec w$ means that $z<w$ and $\mu(z,w)\ne 0$.

\medskip

\def\vp{\varphi}
\def\st{\stackrel}
\def\sc{\scriptstyle}

\medskip

\noindent{\bf 1.2. Cells of $W$.} \quad We refer to [KL] for definition of left cell, right cell and two-sided cell of $W$. For $h,\, h'\in H$ and $x\in W$, write
\begin{alignat*}{2} hC_x&=\sum_{y\in W}a_yC_y,\qquad  C_xh=\sum_{y\in W}b_yC_y,\qquad
  hC_xh'=\sum_{y\in W}c_yC_y,\end{alignat*}
  where $a_y,\, b_y,\, c_y\in \mathcal A.$
  Define $y\ll x$ if $a_y\ne 0$ for some $h\in H$, $y\rl x$ if $b_y\ne 0$ for some $h\in H$, and $y\lrl x$ if $c_y\ne 0$ for some $h,h'\in H$.

  We write $x\el y$ if $x\ll y\ll x$, $x\er y$ if $x\rl y\rl x$, and $x\elr y$ if $x\lrl y\lrl x$. Then $\el,\ \er,\ \elr$ are equivalence relations on $W$. The equivalence classes are called left cells, right cells, and two-sided cells of $W$ respectively.
\medskip

For $w\in W$, set $R(w)=\{s\in S\,|\, ws\le w\}$ and $L(w)=\{s\in S\,|\, sw\le w\}.$ Then we have (see [KL])

\medskip

(a) $R(w)\subset R(u)$ if $u\ll w$ and $L(w)\subset L(u)$ if $u\rl w.$ In particular, $R(w)= R(u)$ if $u\el w$ and $L(w)= L(u)$ if $u\er w.$

\medskip

\noindent{\bf 1.3.   $*$-operations}\ \ The $*$-operation introduced in [KL] and   generalized in [L1] is a useful tool in the theory of cells of Coxeter groups. Let $s,t$ be simple reflections in $S$ and assume that $st$ has order  $m\ge 3$. Let $w\in W$ be such that $sw\ge w,\ tw\ge w$. The $m-1$ elements $sw,\ tsw,\ stsw,\ ...,$ is called a left string (with respect to $\{s,t\})$, and the $m-1$ elements $tw,\ stw,\ tstw,\ ..., $ is also called a left string (with respect to $\{s,t\})$. Similarly we define right strings (with respect to $\{s,t\}$).

\medskip
 Assume that $x$ is in a left (resp. right ) string (with respect to $\{s,t\})$ of length $m-1$ and is the $i$th element of the left (resp. right) string,  define ${}^*x$ (resp. $x^*$) to be the $(m-i)$th element of the string, where $*=\{s,t\}$. The following result is proved in [X2].

 \medskip

 (a) Let $x$ be in $W$ such that $x$ is in a left string with respect to $*=\{s,t\}$ and is also in a right string with respect to $\star=\{s',t'\}$. The ${}^*x$ is in a right string with respect to $\{s',t'\}$ and $x^\star$ is in a left string with respect to $\{s,t\}$. Moreover ${}^*(x^\star)=({}^*x)^\star$. We shall write ${}^*x^\star$
 for ${}^*(x^\star)=({}^*x)^\star$.

 \medskip

Following Lusztig [L1] we set $\tilde\mu(y,w)=\mu(y,w)$ if $y< w$ and $\tilde\mu(y,w)=\mu(w,y)$ if $w<y$. For convenience we also set $\tilde\mu(y,w)=0$ if $y\nless w$ and $w\nless y$. Assume that $x_1,x_2,...,x_{m-1}$  and $y_1,y_2,...,y_{m-1}$ are two left strings with respect to $*=\{s,t\}$. Define
$$a_{ij}=\begin{cases}\tilde\mu(x_i,y_j),\quad &\text{if\ } \{s,t\}\cap L(x_i)=\{s,t\}\cap L(y_j),\\
0,\quad &\text{otherwise}.\end{cases}$$
Lusztig proved the following identities (see subsection 10.4 in [L1]).

\medskip

(b) If $m=3$, then $a_{11}=a_{22}$ and $a_{12}=a_{21}$. In other words, $\tilde\mu(x, y)=\tilde\mu({}^*x, {}^*y)$ for any $x, y$ in any two left strings with respect to $*=\{s,t\}$.

\medskip

(c) If $m=4$, then
$$ a_{11}=a_{33},\ a_{13}=a_{31},\ a_{22}=a_{11}+a_{13},\ a_{12}=a_{21}=a_{23}=a_{32}.$$ In particular, $\tilde\mu(x, y)=\tilde\mu({}^*x, {}^*y)$ for any $x, y$ in any two  left strings with respect to $*=\{s,t\}$.

\medskip

\noindent{\bf 1.4. Lusztig's $a$-function}\quad
For $x,y\in W$, write $C_xC_y=\sum_{z\in
W}h_{x,y,z}C_z,\ h_{x,y,z}\in \mathcal A
 .$  Following Lusztig ([L1]), we define
 $$a(z)={\rm min}\{i\in\bold N\ |\ q^{-\frac i2}h_{x,y,z}\in\mathbb Z[q^{-\frac
 12}]{\rm\ for \ all\ }x,y\in W\}.$$
 The following properties are proved in [L1].

 \medskip

 (a)  We have $a(w)\le
l(w_0)$ for any $w\in W$, where $w_0$ is the longest element in the Weyl group $W_0$.

\medskip

(b) $a(x)\ge a(y)$ if $x\lrl y$. In particular, $a(x)=a(y)$ if $x\elr y$.

\medskip

(c) $x\el y$ (resp. $x\er y,\ x\elr y$) if $a(x)=a(y)$ and $x\ll y$ (resp. $x\lr y,\ x\llr y$).

(d) If $h_{x,y,z}\ne 0$, then $z\lr x$ and $z\ll y$. In particular, $a(z)>a(x)$ if $z\not \er x$, and $a(z)>a(y)$ if $z\not \el y$.

Following Lusztig,
 we define  $\gamma_{x,y,z}$ by $\displaystyle h_{x,y,z}=\gamma_{x,y,z}q^{\frac {a(z)}2}+
 {\rm\ lower\ degree\ terms}.$
Springer showed that $l(z)\ge a(z)$ (see [L2]). Let $\delta(z)$ be
the
 degree of $P_{e,z}$, where $e$ is the neutral element of $W$.
 Then actually one has $l(z)-a(z)-2\delta(z)\ge 0$ (see [L2]). Set
$\cd =\{z\in W\ |\ l(z)-a(z)-2\delta(z)=0\}.$
 The elements of $\cd$ are involutions, called distinguished involutions of
 $(W,S)$ (see [L2]).

\medskip
\def\tt{\tilde T}

For $w\in W$, set $\tilde T_w=q^{-l(w)/2}T_w$. For $x,y\in W$, write
$$\tt_x\tt_y=\sum_{z\in W}f_{x,y,z}\tt_z,\qquad f_{x,y,z}\in\mathcal A=\mathbb Z[q^{\frac12},q^{-\frac12}].$$

(e) If $x,y,w$ are in a two-sided cell of $W$ and $f_{x,y,w}=\lambda q^{\frac{a(w)}2}+$ lower degree terms and as Laurent polynomials in $q^{\frac12}$, deg $f_{x,y,z}\le a(w)$ for all $z\in W$, then
$\gamma_{x,y,w}=\lambda.$

  (f) $x\el y^{-1}$ if and only if $\gamma_{x,y,z}\ne 0$ for some $z\in W$.

 (g) $\gamma_{x,y,z}=\gamma_{y,z^{-1},x^{-1}}=\gamma_{z^{-1}, x,y^{-1}}$.

(h) $\gamma_{x,d,x}=\gamma_{d,x^{-1},x^{-1}}=\gamma_{x^{-1},x,d}=1$ if $x\el d$ and $d$ is a distinguished involution.

 (i)  $\gamma_{x,y,z}=\gamma_{y^{-1},x^{-1},z^{-1}}.$

 (j) If $\omega,\tau\in W$ has length 0, then
 $$\gamma_{\omega x,y\tau,\omega z\tau}=\gamma_{x,y,z},\ \
 \gamma_{ x\omega,\tau y, z}=\gamma_{x,\omega\tau y,z}.$$

 (k) Let $x,y,z\in W$ be such that (1) $x$ is in a left string with respect to $*=\{s,t\}$ and also in a right string with respect to $\#=\{s',t'\}$, (2) $y$ is in a left string with respect to $\#=\{s',t'\}$ and also  in a right string with respect to $\star=\{s'',t''\}$, (3) $z$ is in a left string with respect to $*=\{s,t\}$ and also  in a right string with respect to $\star=\{s'',t''\}$. Then
 $$\gamma_{x,y,z}=\gamma_{{}^*x^\#,{}^\#y^\star,{}^*z^\star}.$$

(l) Let $I$ be a subset of $S$ such that the subgroup $W_I$ of $W$ generated by $I$ is finite. Then the longest element $w_I$ is a distinguished involution.

{ (m) Each left cell (resp. each right cell) of $W$ contains a unique distinguished involution.}

\medskip

Let $d$ is a distinguished involution in $W$.

(n) For any $\omega\in\Omega$, the element $\omega d\omega^{-1}$ is a distinguished involution.

(o) Suppose $s,t\in S$ and $st$ has order 3. Then $d\in D_{L}(s,t)$ if and only if $d\in D_{R}(s,t)$. If $d\in D_{L}(s,t)$, then ${}^*d^*$ is a distinguished involution.

\medskip

{\bf 1.5.} Assume $s, t\in S$ and $st$ has order 4. Let $w, u, v$ be in $W$ such that $l(ststw)=l(w)+4$ and $l(ststv)=l(v)+4$. We have (see [X2 1.6.3])
\begin{itemize}
\item[(a)] $\gamma_{tsw, u, tv}= \gamma_{sw, u, stv},$
\item[(b)] $\gamma_{tsw, u, tsv}= \gamma_{sw, u, sv}+ \gamma_{sw, u, stsv},$
\item[(c)] $\gamma_{tsw, u, tstv}= \gamma_{sw, u, stv},$
\item[(d)]$\gamma_{tstw, u, tv}+\gamma_{tw, u, tv}=\gamma_{stw, u, stv},$
\item[(e)]$\gamma_{tstw, u, tsv}=\gamma_{stw, u, stsv},$
\item[(f)]$\gamma_{tstw, u, tstv}+\gamma_{tw, u, tstv}=\gamma_{stw, u, stv}.$
\end{itemize}

\medskip

Assume $s, t\in S$ and $st$ has order 4. Let $w, u, v$ be in $W$ such that $l(ustst)=l(u)+4$ and $l(vstst)=l(v)+4$. We have (loc.cit.)
\begin{itemize}
\item[(a')] $\gamma_{w, ut, vst}= \gamma_{w, uts, vs},$
\item[(b')] $\gamma_{w, ust, vst}= \gamma_{w, us, vs}+ \gamma_{w, usts, vs},$
\item[(c')] $\gamma_{w, utst, vst}= \gamma_{w, uts, vs},$
\item[(d')]$\gamma_{w, ut, vtst}+\gamma_{w, ut, vt}=\gamma_{w, uts, vts},$
\item[(e')]$\gamma_{w, ust, vtst}=\gamma_{w, usts, vts},$
\item[(f')]$\gamma_{w, utst, vtst}+\gamma_{w, utst, vt}=\gamma_{w, uts, vts}.$
\end{itemize}

\medskip

\def\ll{\underset {L}{\leq}}
\def\rl{\underset {R}{\leq}}

\def\lrl{\underset {LR}{\leq}}
\def\llr{\lrl}
\def\el{\underset {L}{\sim}}
\def\er{\underset {R}{\sim}}
\def\elr{\underset {LR}{\sim}}
\def\ds{\displaystyle\sum}

\def\vp{\varphi}
\def\st{\stackrel}
\def\sc{\scriptstyle}

 \noindent{\bf 1.6. The based ring of a two-sided cell}\quad For each two-sided cell $\Omega$ of $W$, let $J_\Omega$ be the free $\mathbb Z$-module with a basis $t_w,\ w\in \Omega$. Define
  $$t_xt_y=\sum_{z\in \Omega}\gamma_{x,y,z}t_z.$$
  Then $J_\Omega$ is an associative ring with unit $\sum_{d\in\cd\cap \Omega}t_d.$ {We have $J_\Omega=span\{t_w\ |\ w\in\Omega\}=(\sum_{d\in\cd\cap\Omega}t_d)J(\sum_{d\in\cd\cap\Omega}t_d)$. If  $\Gamma\subset\Omega$ is a left cell, and $d_\Gamma$ is the unique distinguished involution in $\Gamma$, then
$$J_{\Gamma\cap\Gamma^{-1}}=t_{d_\Gamma}Jt_{d_\Gamma}\subset J_\Omega$$ is also a based ring, with identity element $t_{d_\Gamma}$.
If $\Gamma, \Theta$ are two left cells, then $t_{d_\Gamma}Jt_{d_{\Theta}}=J_{\Gamma\cap{\Theta}^{-1}}$ is a based $t_{d_\Theta}Jt_{d_\Theta}-t_{d_\Gamma}Jt_{d_\Gamma}$-bimodule.}

  The ring $J=\bigoplus_{\Omega}J_\Omega$ is a based ring with unit $\sum_{d\in\cd}t_d$. Sometimes $J$ is called the asymptotic Hecke algebra since Lusztig established an injective $\ca$-algebra homomorphism (see [L2])
  \begin{alignat*}{2}   \phi: H&\to J\otimes\ca,\\
  C_x&\mapsto\sum_{\st {\st {d\in\cd}{w\in W}}{w\el d}}h_{x,d,w}t_w.\end{alignat*}

  \section{The two-sided cell $c$ of $W$}

  Keep the notations in previous parts. Then  $G=Sp_6(\mathbb C)$, and  $W$ is the extended affine Weyl group attached to $G$ with type $\tilde B_3$. The left cells and two-sided cells are described by J. Du (see [D]). We recall some of his results.

  {\bf 2.1. The Coxeter graph of $W$}. As usual, we number the 4 simple reflections $s_0,\ s_1,\ s_2,\ s_3$ in $W$ so that
  \begin{alignat*}{2} &s_0s_1=s_1s_0,\quad s_0s_3=s_3s_0,\quad s_1s_3=s_3s_1,\\
  &(s_0s_2)^3=(s_1s_2)^3=e,\quad (s_2s_3)^4=e,\end{alignat*}
  where $e$ is the neutral element in $W$. The  Coxeter graph is:

 
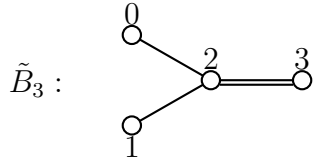
\begin{center}
  \begin{tikzpicture}[scale=.6]
    \draw (-1,0) node[anchor=east]  {$\tilde B_3:$};
    \draw[thick] (2 cm,0) circle (.2 cm) node [above] {$2$};
    \draw[xshift=2 cm,thick] (150:2) circle (.2 cm) node [above] {$0$};
    \draw[xshift=2 cm,thick] (210:2) circle (.2 cm) node [below] {$1$};
    \draw[thick] (4 cm,0) circle (.2 cm) node [above] {$3$};
    \draw[xshift=2 cm,thick] (150:0.2) -- (150:1.8);
    \draw[xshift=2 cm,thick] (210:0.2) -- (210:1.8);
    \draw[thick] (2.2,0) --+ (1.6,0);
    \draw[thick] (2.2,-0.1) --+ (1.6,0);
 \end{tikzpicture}
\end{center}

  There is a unique nontrivial element $\tau $ in $W$ with length 0. We have $\tau^2=e,\ \tau s_0\tau=s_1,\ \tau s_i\tau=s_i$ for $i=2,3.$ Note that $s_1,s_2,s_3$ generate the Weyl group $W_0$ of type $B_3$ and $s_0,s_1,s_2,s_3$ generate an affine Weyl group $W'$ of type $\tilde B_3$. And $W$ is generated by $\tau,\ s_0,s_1,s_2,s_3$.

  \medskip

  {\bf 2.2. The two-sided cell $c$ }\quad According to [D],  $W$  has 8 two-sided cells. We shall use $c$ to denote the two-sided cell of $W$ containing $s_2s_3s_2s_3$, so that $a(c)=4$.

  The two-sided cell $c$ corresponds to the  unipotent class in $Sp_6(\mathbb C)$ with Jordan blocks (2211) (loc.cit.).
  A maximal reductive subgroup of the centralizer of an element in the unipotent class is $F= SL_2(\mathbb C)\times O_2(\mathbb C)\backsimeq SL_2\times (\mathbb Z/2\mathbb Z\ltimes\mathbb C^*)\backsimeq\mathbb Z/2\mathbb Z\ltimes (SL_2\times\mathbb C^*)$, where $\mathbb Z/2\mathbb Z$ acts on $SL_2(\mathbb C)$ trivially and for the nontrivial element $g\in \mathbb Z/2\mathbb Z$ and $z\in\mathbb C^*$ we have $gzg=z^{-1}$.

  \medskip

  {\bf 2.3. Left cells in $c$}\quad  According to  [D, Figure I, Theorem 6.4], $c$ has 18 left cells. We list the 18 left cells and representative elements in the left cells given in [D, Figure I]:

$$\begin{array}{lllllllllll}
&\Gamma_{23}, &2323;& \Gamma_{03},& 23230; &\Gamma_{13},  &23231;
&\Gamma_{02};&232302;\\
&&&&&&&&\\
&\Gamma_{013}, &232301; &\Gamma_{12}, &232312;
 &\Gamma_{01},&2323021;  &\Gamma_2,&2323012;\\
 &&&&&&&&\\
 & {\Gamma'_{01}}, &2323120;
  &\Gamma'_{12}, &23230212;&\Gamma_3,&23230123; &\Gamma'_{02}, &23231202;\\
  &&&&&&&&\\
  &\Gamma'_{13},&232302123; &\Gamma'_{03},&232312023;&\Gamma'_2,&2323021232;
 & {\Gamma''_2},  &2323120232;\\
 &&&&&&&&\\
   &\Gamma_{0}, &23230212320; &\Gamma_{1},&23231202321.&&&& \end{array}$$

 Here  we  simply write $i_1i_2\cdots i_k$ for a reduced expression  $s_{i_1}s_{i_2}\cdots s_{i_k}$ of an element in $W$.

 {\bf 2.4.}  Let $\Gamma$ and $\Gamma'$ be two left cells of $W$. If $\Gamma'=\Gamma^*$ for some $*=\{s,t\}$ (see subsection 1.3 for definition of $*$-operation), then we write $\Gamma\ \overset{\{s, t\}}{\text{------}}\ \Gamma'$. The following result is easy to verify.

$$ \Gamma_{03}\ \overset{\{s_0, s_2\}}{\text{------}}\ \Gamma_{23}\ \overset{\{s_1, s_2\}}{\text{------}}\ \Gamma_{13};\qquad \Gamma_{2}\ \overset{\{s_0, s_2\}}{\text{------}}\ \Gamma_{013}\ \overset{\{s_2, s_3\}}{\text{------}}\ \Gamma_{3};$$

$$\begin{aligned} \Gamma_{02}\ \overset{\{s_1, s_2\}}{\text{------}}\ \Gamma_{01}\ \overset{\{s_0, s_2\}}{\text{------}}\ \Gamma'_{12}\ \overset{\{s_2, s_3\}}{\text{------}}\ &\Gamma'_{2}\ \overset{\{s_0, s_2\}}{\text{------}}\ \Gamma_{0}\\
&\ | {\scriptstyle{\{s_1, s_2\}}}\\
&\Gamma'_{13};\end{aligned}$$

$$\begin{aligned} \Gamma_{12}\ \overset{\{s_0, s_2\}}{\text{------}}\ \Gamma'_{01}\ \overset{\{s_1, s_2\}}{\text{------}}\ \Gamma'_{02}\ \overset{\{s_2, s_3\}}{\text{------}}\ &\Gamma''_{2}\ \overset{\{s_1, s_2\}}{\text{------}}\ \Gamma_{1}\\
&\ | {\scriptstyle{\{s_0, s_2\}}}\\
&\Gamma'_{03}.\end{aligned}$$

 In addition, we have $\Gamma_0=\tau\Gamma_1\tau.$

 \medskip

 {\bf 2.5.} For two left cells $\Gamma$ and $\Gamma'$ of $W$, we say that $\Gamma$ and $\Gamma'$ are {\bf joined} if $\Gamma'$ can be obtained from $\Gamma$ by a sequence of right star operations or  from $\Gamma\tau$ by a sequence of right star operations.

\medskip

{\bf Lemma 2.6.}  Let $\Gamma,\ \Gamma',\ \Theta,\ \Theta',\ \Omega,\ \Omega'$ be left cells of $W$ such that $\Gamma$ and $\Gamma'$ are joined, $\Theta$ and $\Theta'$ are joined and $\Omega$ and $\Omega'$  are joined.
We have:

(a) Through star operations or/and multiplied by $\tau$ we get natural bijections
 $$\Gamma\cap\Theta^{-1}\to \Gamma'\cap\Theta'^{-1},\ \Omega\cap\Gamma^{-1}\to \Omega'\cap\Gamma'^{-1},\  \Omega\cap \Theta^{-1}\to \Omega'\cap \Theta'^{-1}.$$

 (b) Let $x\in\Gamma\cap\Theta^{-1},\ y\in\Omega\cap\Gamma^{-1},\ z\in \Omega\cap \Theta^{-1}$ and $x'\in\Gamma'\cap\Theta'^{-1},\ y'\in\Omega'\cap\Gamma'^{-1},\ z'\in \Omega'\cap \Theta'^{-1}$ be the corresponding elements under the bijections in (a). Then
 $$\gamma_{x,y,z}=\gamma_{x',y',z'}.$$

 {\bf Proof.} It is sufficient to prove assertion (a) and (b)  for $\Gamma'=\Gamma,\ \Gamma^*,\ \Gamma\tau$ and $\Theta'=\Theta,\ \Theta^*,\ \Theta\tau$. According to definition of star operations and $\tau^2=e$, using 1.4(j) and 1.4(k) we see in these cases the assertions are trivial. The lemma is proved.

\medskip

 In subsequent sections, for a reduced expression  $s_{i_1}s_{i_2}\cdots s_{ik}$ of an element in $W$, we often write $i_1i_2\cdots i_k$ instead of the reduced expression.

\section{The based ring of two-sided cell containing $s_2s_3s_2s_3$}

{\bf 3.1.} Keep the notations in section 2.   Let
$$\begin{aligned}
Y_1=\{&\Gamma_{23}, \ \Gamma_{03},\  \Gamma_{13}\}; \qquad gY_1=\{g\Gamma_{23}, \ g\Gamma_{03},\  g\Gamma_{13}\}; \\
Y_2=\{&\Gamma_{013}, \ \Gamma_2, \ \Gamma_3,\}; \qquad gY_2=\{g\Gamma_{013}, \ g\Gamma_2, \ g\Gamma_3,\};\\
Y_3=\{&\Gamma_{02},\ \Gamma_{01},\ \Gamma'_{12},\ \Gamma'_2,\ \Gamma'_{13},\ \Gamma_0,\ \Gamma_{12},\ \Gamma'_{01},\ \Gamma'_{02},\ \Gamma''_2,\ \Gamma'_{03},\ \Gamma_1\};\end{aligned}$$
where $g$ is the non-trivial element in $\mathbb Z/2\mathbb Z$.

 Let $F= SL_2(\mathbb C)\times O_2(\mathbb C)$ act on $Y_3$ trivially and let $F^\circ=SL_2(\mathbb C)\times\mathbb C^*$ act on $Y_1$ and $Y_2$ trivially. Then $F$ acts on the following set algebraically:
 $$Y=Y_1\cup gY_1\cup Y_2\cup gY_2\cup Y_3.$$

Let $Z_{\alpha}=\{\Gamma_\alpha, g\Gamma_\alpha\}$ for any $\Gamma_{\alpha}\in Y_1\cup Y_2$, where $\alpha$ stands for a sequence of digits.

Let  $V(k)$ be an irreducible  rational representation of $SL_2(\mathbb C)$ with highest weight $k$, $U(k)$ be an irreducible  rational representation of $O_2(\mathbb C)$ of dimension 2 with highest weight $k$, $\epsilon$ be the nontrivial irreducible representation of $O_2(\mathbb C)$ of dimension 1, $\eta^k$ be an irreducible  rational representation of $\mathbb C^*$ of dimension 1 with  weight $k$. We often use {\bf 1} for $V(0)$. The representations $V(k)$ and $\eta^k$ are naturally regarded as representations of $F^\circ$. Also the representations $V(k)$, $U(k)$ and $\epsilon$ are naturally regarded as representations of $F$.

We shall denote Irr\,$F$ the isomorphism classes of irreducible rational representations of $F$ and Rep\,$F$ the representation ring of $F$.

 The main result in this section is the following theorem which verifies Lusztig's conjecture on the structure $J_c$.

 \medskip

 {\bf Theorem 3.2.} There exists a bijection
  $$\pi: c\to \text{the set of isomorphism classes of irreducible $F$-vector bundles on}\ Y\times Y,$$
  such that

  (a)  The bijection $\pi$ induces a based ring isomorphism
  $$\pi: J_c\to K_{F}(Y\times Y),\ \ t_x\mapsto \pi(x).$$

  (b) $\pi(x^{-1})_{(a,b)}=\pi(x)_{(b,a)}^*$ is the dual representation of $\pi(x)_{(b,a)}.$

  \medskip

  {\bf Remark.} This result actually has already been proved in Theorem 4 in [BO] conceptually. Here we will construct the bijection explicitly so that the isomorphism can be used to compute certain irreducible representations of affine Hecke algebras of type $\tilde B_3$.

  \medskip

 The rest of this section is devoted to proving the theorem above. The following convention will be frequently used.

 {\bf Convention:} {\it The symbol $\Box$ will stand for any element in the two-sided ideal $H^{<2323}$ of $H$ spanned by all $C_w$ with $a(w)> 4$.} As a result, for any $h\in H$, we have $h\Box=\Box h=\Box,\ \Box+\Box=\Box$.

 \medskip

  {\bf 3.3. }   In this subsection we consider $\Gamma_{02}\cap\Gamma_{02}^{-1}$ and its based ring $J_{\Gamma_{02}\cap\Gamma_{02}^{-1}}$.

   Let $*=\{s_0, s_2\}$, $\sharp=\{s_1, s_2\}$ and $\star=\{s_2, s_3\}$.  According to [D, Theorem 6.4], the set $\Gamma_{02}\cap\Gamma^{-1}_{02}$ consists of the following elements:
 $$\begin{aligned}
 d_{02}&=s_2s_0s_3s_2s_3s_2s_0s_2,\qquad x_0:=s_2s_0s_2s_3s_2s_0;\\
 x_k&=s_2s_0s_3s_2(s_3s_2s_1s_2\tau)^k, \qquad k\ge 1\\
y_k&=s_2s_0s_3s_2s_3s_2s_0s_2(s_1s_2s_3s_2s_1\tau)^k=d_{02}(s_1s_2s_3s_2s_1\tau)^k, \quad k\ge 0 \\
 \dot y_k&={}^\star y_k=y_k^\star=x_0(s_1s_2s_3s_2s_1\tau)^k, \quad  k\geq 0.\\
 z_{i,k}&:=\begin{cases}{x_1(s_1s_2s_3s_2s_1\tau)^k,}&{i=1,\ k\ge 1;}\\
 {\tau s_2(^{\sharp*}(s_3z_{i-1,k})),}&{i\geq 2,\ k\ge 1.}
\end{cases}
\end{aligned}$$

\def\Irr{\text{Irr}}
\def\Rep{\text{Rep}}

  \medskip

  {\bf Proposition 3.3A.} Keep the notations in subsection 3.1. That is,
 $c$ stands for the two-sided cell of $W$ (the extended affine Weyl group attached to $Sp_6(\mathbb C)$) containing $s_2s_3s_2s_3$, and $\Gamma_{02}$ is the left cell in $c$ containing $s_2s_3s_2s_3s_0s_2$. Then the bijection
\begin{alignat*}{7}
\pi: \Gamma_{02}\cap \Gamma^{-1}_{02}&\longrightarrow \Irr F, &&&&&\\
d_{02}&\longmapsto V(0),&\ \  x_0&\longmapsto \epsilon,&&&\\
x_k&\longmapsto U(k),&\ \ y_k&\longmapsto  V(k),&\ \ k\ge 1;&\\
 \dot y_k&\longmapsto \epsilon\otimes V(k),&
\ \ z_{{i},k}&\longmapsto U(i)\otimes V(k), &\ \ k\ge 1,\ i\ge 1;\end{alignat*}
induces a based ring isomorphism
$$\pi: J_{\Gamma_{02}\cap \Gamma^{-1}_{02}}\to \Rep\,F,\quad t_{w}\mapsto \pi(w).$$

\medskip

Clearly $\pi$ is bijective. To prove the proposition we need to verify $\pi(t_wt_u)=\pi(t_w)\pi(t_u)$. \textbf{The product $\pi(t_w)\pi(t_u)$ is easy to calculate, so the main job is to compute $t_wt_u$.} We prove the proposition by establishing a series of lemmas. Note that for any $x$ in $\Gamma_{02}\cap \Gamma^{-1}_{02}$, we have $x=x^{-1}$. Using 1.4 (i) we get

    \medskip

    {\bf Lemma 3.3B.}  $t_wt_u=t_ut_w$ for any $w, u\in\Gamma_{02}\cap \Gamma^{-1}_{02}$.

     \medskip

{\bf Lemma 3.3C.} (1) The element $d_{02}$ is a distinguished involution and $x_0=d_{02}^\star={}^\star d_{02}$ (recall that $\star=\{s_2,s_3\}$).

(2) Let $w=d_{02}$ or $x_0$. Then we have $\pi(t_wt_u)=\pi(t_w)\pi(t_u)$ for any $u$ in $\Gamma_{02}^{-1}\cap\Gamma_{02}$.

\begin{proof} As we explained, the main job is to compute $t_wt_u$.

For simplicity we write $d$ for $d_{02}=s_2s_0s_3s_2s_3s_2s_0s_2.$ As usual, $e$ denotes the neutral element of $W$. By a direct computation we get  $P_{e,d}=P_{s_2s_0s_2,d}=1+q^2$. So the  degree $\delta(d)$ of $P_{e,d}$ is 2. Then  $l(d)-a(d)-2\delta(d)=8-4-4=0$. Hence $d$ is a distinguished involution (see subsection 1.4). Apply 1.4 (h) we see that $t_dt_u=t_u$ for any $u\in\Gamma_{02}\cap\Gamma_{02}^{-1}$.

Since $d=s_2s_3 s_2s_0s_2s_3s_2s_0=s_0s_2s_3s_2s_0s_2s_3s_2$, hence, for $\star=\{s_2,s_3\}$, we have ${}^\star d=d^\star=s_2s_0s_2s_3s_2s_0=x_0.$  Note that ${}^\star x_k=x_k$ for $k\ge 1$, ${}^\star y_k=\dot y_k$ for $k\ge 0$, and ${}^\star z_{i,k}=z_{i,k}$ for $i, k\ge 1$. Applying 1.4 (h) and 1.4 (k), from $t_dt_u=t_u$ we can see easily  $\pi(t_{x_0}t_u)=\pi(t_{x_0})\pi(t_u)$ for any $u$ in $\Gamma_{02}^{-1}\cap\Gamma_{02}$.
\end{proof}

\medskip

{\bf Lemma 3.3D.} Let  $l$ and $k$ be  positive integers. Then $\pi(t_{x_l}t_{x_k})=\pi(t_{x_l})\pi(t_{x_k})$.

\medskip

\begin{proof} It suffices to prove  (a) $t_{x_k}t_{x_k}=t_{x_{2k}}+t_{d_{02}}+t_{x_0},$
(b) $t_{x_l}t_{x_k}=t_{x_{k+l}}+t_{x_{| k-l |}}$ if $k\ne l$.

Let $\xi=q^{\frac12}+q^{-\frac12}$. By a simple computation we get
\begin{alignat}{2}C_{x_1}&=C_\tau (C_{s_2}C_{s_1}-1)(C_{s_3}C_{s_2}-2)C_{s_3}C_{s_2s_0s_2}.\end{alignat}

Since $C_{s_2s_0s_2}C_{x_k}=(\xi^3-\xi)C_{x_k}$, we get
\begin{align}
C_{x_1}C_{x_k}=(\xi^3-\xi)C_\tau (C_{s_2}C_{s_1}-1)(C_{s_3}C_{s_2}-2)C_{s_3}C_{x_k}.\end{align}

 By a direct computation we get for any $k\ge 2$,
\begin{equation}\begin{aligned}
C_{x_1}C_{x_k}
=&(\xi^4-\xi^2)(C_{x_{k+1}}+C_{x_{k-1}})+\Box\in (\xi^4-\xi^2)(C_{x_{k+1}}+C_{x_{k-1}})+H^{<2323},
\end{aligned}\end{equation}
where  $\Box$ is an element in $H^{<2323}$ (see the Convention below Theorem 3.2).

\medskip

 The formula above implies the following identity.
 \begin{align}t_{x_1}t_{x_k}=t_{x_{k-1}}+t_{x_{k+1}},\qquad  k\ge 2.\end{align}

Now we deal with the case $k=1$. By formula  (2) we have
\begin{align}
C_{x_1}C_{x_1}=(\xi^3-\xi)C_\tau (C_{s_2}C_{s_1}-1)(C_{s_3}C_{s_2}-2)C_{s_3}C_{x_1}.
\end{align}

Again by a direct computation we get
$$\begin{aligned}
C_{x_1}C_{x_1}=&(\xi^4-\xi^2)(C_{x_2}+C_{d_{02}}+C_{x_0})+\Box.
\end{aligned}$$

Therefore we have  \begin{align}t_{x_1}t_{x_1}=t_{x_2}+t_{d_{02}}+t_{x_0}.\end{align}

\medskip

For $1\le l\le k-1$, using formulas (4), (6), induction on $l$ and Lemma 3.3C we get
$$t_{x_l}t_{x_k}=t_{x_{k+l}}+t_{x_{k-l}}.$$
Noting that $t_{x_k}t_{x_l}=t_{x_l}t_{x_k}$ by Lemma 3.3B, hence for any different positive integers $k$ and $l$, (b) holds.

Using induction on $k$, above formula, formula (6) and Lemma 3.3C, we see
$$t_{x_k}t_{x_k}=t_{x_{2k}}+t_{d_{02}}+t_{x_0}.$$

The proof is completed.
\end{proof}

\medskip

{\bf Lemma 3.3E.} Let $x=x_i, y=y_k$ or $\dot y_k$ with $i, k$ being positive integers. Then
$\pi(t_xt_y)=\pi(t_x)\pi(t_y)$.

\medskip

\begin{proof} Recall $\star=\{s_2,s_3\}$. Note  $x_k=x_k^\star$ and $\dot y_k={}^\star y_k$. By 1.4(k), it suffices to show
$t_{x_i}t_{y_k}=t_{z_{i,k}}$.

 First we  prove  (a) $t_{x_1}t_{y_k}=t_{z_{1,k}},$  (b) $t_{s_2x_1}t_{y_k}=t_{s_2z_{1,k}}$.

Using 1.4(k) we see that (b) is equivalent to (b') $t_{s_3x_1}t_ {y_k}=t_{s_3z_{1,k}}$ since $s_3x_1={}^\star(s_2x_1)$ and $s_3z_{1,k}={}^\star(s_2z_{1,k})$.

As before $\xi=q^{\frac12}+q^{-\frac12}$. For $k\geq 1$, by formula (1) and $C_{s_2s_0s_2}C_{y_k}=(\xi^3-\xi)C_{y_k}$, we get,
\begin{align}
C_{x_1}C_{y_k}=(\xi^3-\xi)C_\tau(C_{s_2}C_{s_1}-1)(C_{s_3}C_{s_2}-2)C_{s_3}C_{y_k}.
\end{align}

\medskip

By a direct computation we get

\begin{align*}
C_{x_1}C_{y_k}=&(\xi^4-\xi^2)C_{z_{1,k}}+\Box,
\end{align*}
where  $\Box$ is an element in $H^{<2323}$ (see the Convention below Theorem 3.2).

\medskip

So (a) is true.

Let $z\in W$ be such that $\gamma_{s_2x_1, y_k, z}\ne 0$. By 1.4(f) and 1.4(g) we have $z\er s_2x_1$. Since $s_2x_1$ is in a left string with respect to $\{s_2, s_3\}$,  so is $z$. Then using  1.5(a), 1.5(b),1.5(f) and (a) we obtain (b).

\medskip

 Now  assume that $i\ge 2$ and $t_{x_l}t_{y_k}=t_{z_{l,k}}$ for any $l\le i-1$. Similarly using 1.5(a), 1.5(b), 1.5(f) and 1.4(k) we obtain
\begin{align}&t_{s_2x_l}t_ {y_k}=t_{s_2z_{l,k}},\\
&t_{s_3x_l}t_{y_k}=t_{s_3z_{l,k}}.\end{align}

Note $s_3x_{i-1}={}^{*\sharp}(\tau s_2x_i)$. For any $v\in W$, using 1.4(j), 1.4(k) and formula (9), we  see that
$$
\gamma_{s_2x_i, y_k, v}
=\gamma_{{}^{*\sharp}(\tau s_2x_i), y_k, {}^{*\sharp}(\tau v)}
=\gamma_{s_3x_{i-1}, y_k, {}^{*\sharp}(\tau v)}\ne 0
$$
if and only if $v=s_2z_{i,k}$, and $\gamma_{s_2x_i, y_k, s_2z_{i,k}}=1$.

Therefore \begin{align}t_{s_2x_i}t_{y_k}=t_{s_2z_{i,k}}.\end{align}

Let $w\in W$ be such that $\gamma_{x_i, y_k, w}\ne 0$. By 1.4(f) and 1.4(g), $w\er x_i$.  Since $x_i$ is in a left string with respect to $\{s_2,s_3\}$, so is $w$.
Then by 1.5(b), 1.5(c) and formula (10), we get $w=z_{i,k}$ and
$\gamma_{x_i, y_k, z_{i,k}}=1$, so
$t_{x_i}t_{y_k}=t_{z_{i,k}}$. By induction, we get $t_{x_i}t_{y_k}=t_{z_{i,k}}$ for any $i, k\ge 1$.

\medskip

The proof is completed.
\end{proof}

\medskip

{\bf Corollary 3.3F.} Let $i, j, k$ be positive integers. Then
$\pi(t_{x_{j}}t_{z_{i, k}})=\pi(t_{x_{j}})\pi(t_{z_{i, k}})$.
\begin{proof}
It follows Lemma 3.3D and Lemma 3.3E.
\end{proof}

\medskip

{\bf Lemma 3.3G.}  For  any $x, y\in\{\ y_k,\dot y_k\ |\  k\ge 1\}$, we have
$\pi(t_xt_y)=\pi(t_x)\pi(t_y)$.
\begin{proof}
We need to compute the product $t_xt_y$. Note $y_0=d_{02}, \dot y_0=x_0$ and $\dot y_i={}^\star y_i=y_i^\star$ for $\star=\{s_2,s_3\}$, for any $i\ge0$. By 1.4(k) it suffices to prove
$$t_{y_l}t_{y_k}=\sum\limits_{\substack{0\leq i\leq \min \{k, l\}}}t_{y_{k+l-2i}}.$$

We first show the equality above is true for $l=1$, that is
\begin{align}t_{y_1}t_{y_k}=t_{y_{k-1}}+t_{y_{k+1}},\qquad k\ge 1.\end{align}

As before, let $\xi=q^{\frac12}+q^{-\frac12}$ and  set $d:=d_{02}$.
Let $*=\{s_0, s_2\}$, $\sharp=\{s_1, s_2\}$ and $u=\tau s_1 y_1=s_2s_3s_2s_1d={}^\star(^*(^\sharp d))).$
Write
\begin{equation} C_uC_{y_k}=\sum_{z\in W}h_{u,y_k,z}C_z.\end{equation}
According to 1.4(k), 1.4(d) and 1.4(c)  we see

\medskip

(i) If $z\not \er u$, then $a(z)>a(u)=4$.

\medskip

(ii) If $z\er u$, then $z\el y_k$. Moreover, in this case, $\deg h_{u,y_k,z}=4$ (as a Laurent polynomial in $q^{\frac12}$) if and only if
$z={}^\star(^*(^\sharp y_k)))=s_2s_3s_2s_1y_k.$ Furthermore, $\gamma_{u,y_k, {}^\star(^*(^\sharp y_k)))}=1.$

\medskip

Hence

(iii) If $z\er u$ and $z\ne s_2s_3s_2s_1y_k$, then $\deg h_{u,y_k,z}<4$ and $s_1z\ge z$.

\medskip

 Therefore, we get \begin{equation} C_uC_{y_k}=hC_{s_2s_3s_2s_1y_k}+\sum_{\st {\st {s_2s_3s_2s_1y_k\ne z\in W}{z\er u, s_1z>z}}{\deg h_{u,y_k,z}<4}}h_{u,y_k,z}C_z+\Box,\end{equation}
where $ h=(q^{\frac12})^4+$ lower degree terms (as a Laurent polynomial in $q^{\frac12}$).

Note that $u=s_2s_3s_2s_1y_0.$  Now we compute $C_{s_1}C_{s_2s_3s_2s_1y_k}$.

By { 1.1(a)},  we have
\begin{equation}{ C_\tau}C_{s_1}C_{s_2s_3s_2s_1y_k}={ C_{y_{k+1}}}+\sum\limits_{\substack{y\prec s_2s_3s_2s_1y_k\\ s_1y<y}}\mu(y, s_2s_3s_2s_1y_k)C_{\tau y}.\end{equation} Note that ${L}(s_2s_3s_2s_1y_k)=\{s_2\}$.

By a careful analysis, we get
\begin{equation}{ C_\tau}C_{s_1}C_{s_2s_3s_2s_1y_k}={ C_{y_{k+1}}+C_{\tau s_3s_2s_1y_k}+C_{y_{k-1}}}+\Box\end{equation}

Similarly we get
\begin{equation} { C_\tau}C_{s_1}C_u={ C_{y_1}+C_{\tau s_3s_2s_1y_0}}+\Box.\end{equation}

Since { $\tau s_3s_2s_1y_i=\tau{}^{\sharp\star*\sharp}y_i$ for any $i\ge 0$}, by Lemma 3.3C and 1.4(k) we get
\begin{equation} C_{\tau s_3s_2s_1y_0}C_{y_k}=fC_{\tau s_3s_2s_1y_k}+\sum_{\st {\st {s_3s_2s_1y_k\ne z\in W}{z\er \tau s_3s_2s_1y_0}}{\deg f_{\tau s_3s_2s_1y_0,y_k,z}<4}}f_{\tau s_3s_2s_1y_0,y_k,z}C_z+\Box,\end{equation}
where $ f=(q^{\frac12})^4+$ lower degree terms (as a Laurent polynomial in $q^{\frac12}$).

\medskip

Combining formulas (15)-(17) and (13), we get
$$\begin{aligned}C_{y_1}C_{y_k}
=&h(C_{y_{k-1}}+C_{y_{k+1}})+f'C_{\tau s_3s_2s_1y_k}\\
&+(\sum_{\st {\st {s_2s_3s_2s_1y_k\ne z\in W}{z\er u, s_1z>z}}{\deg h_{u,y_k,z}<4}}h_{u,y_k,z}{ C_\tau}C_{s_1}C_z-\sum_{\st {\st {s_3s_2s_1y_k\ne z\in W}{z\er \tau s_3s_2s_1y_0}}{\deg f_{\tau s_3s_2s_1y_0,y_k,z}<4}}f_{\tau s_3s_2s_1y_0,y_k,z}C_z)+\Box,\end{aligned}$$
where $\deg f'<4$ (as a Laurent polynomial in $q^{\frac12}$).

Therefore formula (11) is true. For $l\ge 2$, we have $t_{y_l}=t_{y_1}t_{y_{l-1}}-t_{y_{l-2}}$. Using induction on $l$ and formula (11) we then can compute $t_{y_l}t_{y_k}$ for $l\ge 2$ through the following identity
$$ t_{y_l}t_{y_k}=(t_{y_1}t_{y_{l-1}}-t_{y_{l-2}})t_{y_k}.$$

The proof is completed.
\end{proof}

\medskip

 \medskip

{\bf Corollary 3.3H.} { For $x\in\{y_i, z_{i, k}\ |\ i, k\ge 1\}, y\in\{z_{i, k}\ |\ i, k\ge 1\}$,  we have $\pi(t_xt_y)=\pi(t_x)\pi(t_y)$.}
\begin{proof}
Set $z_{i,0}=x_i$. Then by Lemma 3.3C, Lemma 3.3E, Lemma 3.3G and 1.4(i), we get for any $i, j, k, l\ge 1$,
\begin{itemize}

\item[(a)] $ t_{y_{j}}  t_{z_{i, k}} = t_{\dot y_{j}}  t_{z_{i,k}} =\sum\limits_{\substack{0\leq l\leq min\{j, k\}}} t_{z_{i, j+k-2l}} $;

\item[(b) ]
$ t_{z_{j, l}}  t_{z_{i, k}} =\begin{cases}{\sum\limits_{\substack{0\leq p\leq min\{l, k\}}} t_{z_{2i, l+k-2p}} + t_{y_{l+k-2p}} + t_{\dot y_{l+k-2p}}, }&{j=i}\\
{ \sum\limits_{\substack{0\leq p\leq min\{l, k\}}} t_{z_{i+j,  l+k-2p}} + t_{z_{|i-j|, l+k-2p}}, }&{j\ne i}.\end{cases}$
\end{itemize}

The corollary follows.
\end{proof}
\
Now Proposition 3.3A results from 3.3B-3.3H.

 \medskip

{\bf Corollary 3.3I.}
  For any left cells $\Gamma,\ \Theta$ in $Y_3$ we have a bijection
  $$\Gamma \cap\Theta^{-1}\to\text{Irr}F.$$

 \begin{proof}
It follows from Lemma 2.6 (a), subsection 2.4 and Proposition 3.3A.
\end{proof}

\medskip

 {\bf 3.4.}  In this subsection we consider $\Gamma_{03}\cap\Gamma_{03}^{-1}$ and its based ring $J_{\Gamma_{03}\cap\Gamma_{03}^{-1}}$.

   \medskip

     According to  [D, {Theorem 6.4}], the set $\Gamma_{03}\cap\Gamma^{-1}_{03}$ consists of the elements below:
 $$\begin{aligned}
&\hat x_k:=s_2x_ks_3=s_0s_3s_2s_3(s_2s_1s_2s_3\tau)^k, \hat x_0:=s_2x_0s_3=s_0s_2s_3s_2s_3s_0=s_2d_{02}s_2;\\
&\hat x'_k:=\hat x_k^\star=s_2x_{k}s_2=s_0s_3(s_2s_3s_2s_1\tau)^{k};\\
&\tilde x_k:={}^\star\hat x_k=s_3x_{k}s_3=s_3s_2s_0s_3s_2s_3(s_2s_1s_2s_3\tau)^{k},
\tilde x_0:={}^\star\hat x_0=\hat x_0^\star=s_3x_0s_3=s_0s_2s_3s_2s_3s_0s_2s_3;\\
&\tilde x'_k:={}^\star\hat x_k^\star=s_3x_ks_2=s_3s_2s_0s_3(s_2s_3s_2s_1\tau)^k;\\
&\hat  y_k:=s_2y_ks_2=s_0s_2s_3s_2s_3s_0(s_1s_2s_3s_2s_1\tau)^k;\\
&\hat  y'_k:={}^\star\hat  y_k=\hat  y_k^\star=s_2y_ks_3=s_3y_ks_2=s_0s_2s_3s_2s_3s_0s_2s_3(s_1s_2s_3s_2s_1\tau)^k;\\
&\hat z_{i, k}:=\begin{cases}{s_2z_{1, k}s_3=\hat x_1(s_1s_2s_3s_2s_1\tau)^k,}&{i=1}\\
{\tau ^{\sharp*\star}(\hat z_{i-1, k}),}&{i\geq 2}
\end{cases};\\
&{}^\star \hat z_{i, k},\quad \hat z_{i, k}^\star,\quad {}^\star\hat z_{i, k}^\star,\end{aligned}$$
where $*=\{s_0, s_2\}$, $\sharp=\{s_1, s_2\}$ and $\star=\{s_2, s_3\},\  k\geq 1$.

\medskip
\medskip

  {\bf Proposition 3.4A.} Keep the notations in subsection  3.1. Let   $Z_{03}=\{\Gamma_{03}, g\Gamma_{03}\}$ and $k\ge 1$.
 Then the bijection
$$\pi: \Gamma_{03}\cap \Gamma^{-1}_{03}\longrightarrow \{\text{iso. classes of irr. $F$-vector bundles on}\ Z_{03}\times Z_{03}\}\subset M_{2\times 2}(R_{F^\circ}),$$

\medskip

$
 \hat x_0 \longmapsto   \left(\begin{array}{cc}
\mathbf 1 & 0 \\
0 & \mathbf 1
\end{array}\right),\quad\hat x_k \longmapsto
\left(\begin{array}{cc}
\eta^{-k} & 0 \\
0 & \eta^{k}
\end{array}\right),\quad
\hat x'_k \longmapsto\left(\begin{array}{cc}
0 & \eta^{-k} \\
\eta^{k} & 0
\end{array}\right),$

\medskip

$\tilde x_0 \longmapsto  \left(\begin{array}{cc}
0 & \mathbf 1 \\
\mathbf 1 & 0
\end{array}\right),\quad
\tilde x_k \longmapsto
\left(\begin{array}{cc}
0 & \eta^{k} \\
\eta^{-k} & 0
\end{array}\right),\quad
\tilde x'_k \longmapsto \left(\begin{array}{cc}
\eta^{k} & 0 \\
0 & \eta^{-k}
\end{array}\right),$

\medskip

$\hat y_k \longmapsto  \left(\begin{array}{cc}
V(k) & 0 \\
0 & V(k)
\end{array}\right),\quad
\hat y'_k\longmapsto \left(\begin{array}{cc}
0 & V(k) \\
V(k) & 0
\end{array}\right),$

\medskip

$\hat z_{i,k}  \longmapsto   \left(\begin{array}{cc}
\eta^{-i}\otimes V(k) & 0 \\
0 & \eta^i\otimes V(k)
\end{array}\right),\quad
{}^\star(\hat z_{i,k}) \longmapsto  \left(\begin{array}{cc}
0 & \eta^i\otimes V(k)\\
\eta^{-i}\otimes V(k) & 0
\end{array}\right),$

\medskip

$(\hat z_{i,k})^\star  \longmapsto   \left(\begin{array}{cc}
0& \eta^{-i}\otimes V(k)  \\
\eta^i\otimes V(k) & 0
\end{array}\right),\quad
{}^\star(\hat z_{i,k})^\star \longmapsto  \left(\begin{array}{cc}
\eta^i\otimes V(k) & 0\\
0& \eta^{-i}\otimes V(k)
\end{array}\right),$

\medskip

\noindent  induces a based ring isomorphism
$$\pi: J_{\Gamma_{03}\cap \Gamma^{-1}_{03}}\to  K_{F}(Z_{03}\times Z_{03}),$$

\medskip

Clearly $\pi$ is bijective. To prove the proposition we need to verify $\pi(t_wt_u)=\pi(t_w)\pi(t_u)$. \textbf{The product $\pi(t_w)\pi(t_u)$ is easy to calculate, so the main job is to compute $t_wt_u$.} We prove the proposition by establishing a series of lemmas.

\medskip

{\bf{Lemma 3.4B.}} For any $x\in \{\hat x_0, \tilde x_0\},  y\in\Gamma_{03}\cap\Gamma_{03}^{-1}$,
we have $\pi(t_xt_y)=\pi(t_x)\pi(t_y),$
$ \pi(t_yt_x)=\pi(t_y)\pi(t_x)$.

\begin{proof}
It suffices to show

 (a) $t_{\hat x_0}t_y=t_yt_{\hat x_0}=t_y$,

(b) $t_{\tilde x_0}t_y=t_{{}^\star y}, $

(c)  $t_yt_{\tilde x_0}=t_{y^\star}$.

Since $\hat x_0={}^*(s_2s_3s_2s_3)^*$ for $*=\{s_0, s_2\}$, by 1.4(l) and 1.4(o), we know that $\hat x_0$ is a distinguished involution. Thus  assertion (a) is true and the parts (b) and (c) follow from the fact $\tilde x_0={}^\star \hat x_0=\hat x_0^\star$ and 1.4(k).
\end{proof}

\medskip

{\bf{Lemma 3.4C.}}  For any $x, y\in\{\hat x_k, \hat x'_l, \tilde x_k, \tilde x'_l\ |\  k\ge 0,\ l\ge 1\},$ we have $\pi(t_xt_y)=\pi(t_x)\pi(t_y)$.
\begin{proof}
As we explained before, the main job is to compute the product $t_xt_y$.

Let $x=\hat x_k, y=\hat x_l$,  $k, l\ge 0$. Then in this case it suffices to prove \begin{align}t_{\hat x_k}t_{\hat x_l}=t_{\hat x_{k+l}}.\end{align}

We use induction on $k$ to prove (18). Let $*=\{s_0, s_2\}, \sharp=\{s_1,s_2\}, \star=\{s_2, s_3\}$.

When $k=0$ or $l=0$, formula (18) is a special case of Lemma 3.4B (a).    Let $k=1$. Assume $\gamma_{\hat x_1, \hat x_l, p}\ne0$ for some $p\in\Gamma^{-1}_{03}\cap\Gamma_{03}$. Note $\hat x_k={}^{*\sharp\star}(\tau \hat x_{k-1})$ for any $k\ge 1$. Then by 1.4(j), 1.4(k) and Lemma 3.4B , we have
$$\gamma_{\hat x_1, \hat x_l, p}=\gamma_{{}^{*\sharp\star}(\tau\hat x_{0}), \hat x_l, p}
=\gamma_{\hat x_0, \hat x_l, \tau{}^{\star\sharp*}p}\ne 0$$
if and only if $p={}^{*\sharp\star}(\tau \hat x_{l})=\hat x_{l+1}$. Hence $\gamma_{\hat x_1, \hat x_l, \hat x_{l+1}}=1$ and $\gamma_{\hat x_1, \hat x_l, z}=0$ if $z\ne \hat x_{l+1}$. Thus (18) is valid when $k=1$.

For $k\ge 2$, since $t_{\hat x_k}t_{\hat x_l}=t_{\hat x_1}t_{\hat x_{k-1}}t_{\hat x_l}$, then by induction on $k$ we  can deduce (18) easily.

\medskip

Set $t_{\hat x'_0}=t_{\tilde x'_0}=t_{\hat x_{i}}=t_{\tilde x_{i}}=t_{\hat x'_{i}}=t_{\tilde x'_{i}}=0$ if $i<0$. Similarly, we can prove

(a)  $ t_{\hat x_k}  t_{\hat x'_l} = t_{\hat x'_{k+l}} ,$ $ t_{\hat x_k}  t_{\tilde x_l} = t_{\tilde x_{l-k}} + t_{\hat x'_{k-l}} ,$
$ t_{\hat x_k}  t_{\tilde x'_l} = t_{\tilde x'_{l-k}} + t_{\hat x_{k-l}} ;$

Note $\hat x'_k=\hat x_k^\star,\ \tilde x_k={}^\star\hat x_k,\ \tilde x'_k={}^\star\hat x_k^\star$, the following identities follow from (18) and (a)  using 1.4(k) .

(b)  $ t_{\hat x'_k}  t_{\hat x_l} = t_{\tilde x_{l-k}} +
 t_{\hat x'_{k-l}} $, $ t_{\hat x'_k}  t_{\hat x'_l} = t_{\tilde x'_{l-k}} +
t_{\hat x_{k-l}} $, $ t_{\hat x'_k}  t_{\tilde x_l} = t_{\hat x_{k+l}} $, $ t_{\hat x'_k}  t_{\tilde x'_l} = t_{\hat x'_{k+l}} ;$

(c)  $ t_{\tilde x_k}  t_{\hat x_l} = t_{\tilde x_{k+l}} $, $ t_{\tilde x_k}  t_{\hat x'_l} = t_{\tilde x'_{k+l}} $, $ t_{\tilde x_k}  t_{\tilde x_l} = t_{\hat x_{l-k}} +
 t_{\tilde x'_{k-l}} $, $ t_{\tilde x_k}  t_{\tilde x'_l} = t_{\hat x'_{l-k}} +
t_{\tilde x_{k-l}} $;

(d)   $ t_{\tilde x'_k}  t_{\hat x_l} = t_{\hat x_{l-k}} +
t_{\tilde x'_{k-l}} $, $ t_{\tilde x'_k}  t_{\hat x'_l} = t_{\hat x'_{l-k}} +
t_{\tilde x_{k-l}} $, $ t_{\tilde x'_k}  t_{\tilde x_l} = t_{\tilde x_{k+l}} $, $ t_{\tilde x'_k}  t_{\tilde x'_l} = t_{\tilde x'_{k+l}} $.

\medskip

The proof is completed.
\end{proof}

\medskip

{\bf Lemma 3.4D.}  For any $x\in\{\hat x_i, \hat x'_i, \tilde x_i, \tilde x'_i\ |\ i\ge 1\}$ and any $ y\in\{\hat y_k, \hat y'_k\ |\ k\ge 1\}$, we have $\pi(t_xt_y)=\pi(t_x)\pi(t_y)$.

\begin{proof}
 Note for any $i, k\geq 1$, $\hat x'_i=\hat x_i^\star,\ \tilde x_i={}^\star\hat x_i,\ \tilde x'_i={}^\star\hat x_i^\star, \hat y'_k={}^\star(\hat y_k)=(\hat y_k)^\star$. Using 1.4(k) we see that it suffices to prove $t_{\hat x_i} t_{\hat y_k}= t_{\hat z_{i, k}}$.

 Since $\hat x_{i}=\tau{}^{\sharp *\star}(\hat x_{i-1})$, $\hat z_{i,k}=\tau{}^{\sharp *\star}(\hat z_{i-1,k})$ for any $i\geq 2$, using 1.4(j), 1.4(k) and induction on $k$, we only need to show  \begin{align}t_{\hat x_1}t_{\hat y_k}=t_{\hat z_{1,k}}.\end{align}

 Note $\hat x_{1}=(\hat x_{0})^{\star*\sharp}\tau, \hat y_k=\tau{}^{\sharp*\star}(\hat z_{1,k})$. Using 1.4(j) and 1.4(k), we get for any $p\in W$,
$$\gamma_{{\hat x_1}, {\hat y_k}, p}
=\gamma_{(\hat x_0)^{\star*\sharp}\tau, \tau{}^{\sharp*\star}(\hat z_{1,k}), p}
=\gamma_{\hat x_0, \hat z_{1,k}, p}.$$
By Lemma 3.4B(a), we get $\gamma_{{\hat x_1}, {\hat y_k}, p}\ne 0$ if and only if $p=\hat z_{1,k}$ and $\gamma_{{\hat x_1}, {\hat y_k}, \hat z_{1,k}}=1$. Therefore formula (19) is true.

\medskip

The proof is completed.
\end{proof}

\medskip

{\bf Corollary 3.4E.} For any $x\in\{\hat x_i, \hat x'_i, \tilde x_i, \tilde x'_i\ |\ i\ge 1\}$ and any
$y\in\{\hat z_{j, k}, {}^\star(\hat z_{j, k}), (\hat z_{j, k})^\star,$
$ {}^\star(\hat z_{j, k})^\star\ |\ j, k\ge 1\}$, we have $\pi(t_xt_y)=\pi(t_x)\pi(t_y)$.
\begin{proof}
The result follows from Lemma 3.4C, Lemma 3.4D and 1.4(k).\end{proof}

\medskip

{\bf Lemma 3.4F.}
For  any $x, y\in\{\hat y_k, \hat y'_k \ |\ k\ge 1\}$, we have $\pi(t_xt_y)=\pi(t_x)\pi(t_y)$.
\begin{proof}
Let $\hat y_0=\hat x_0$, $\hat y'_0=\tilde x_0$.

 Note that $\hat y'_k={}^\star \hat y_k=\hat y_k^\star$ for any $k\ge 1$, where $\star=\{s_2,s_3\}$.

Using 1.4(k) we see that it suffices to compute $t_{\hat y_k}t_{\hat y_l}$ with $k, l\ge 1$. We start with the case $k=1$ and prove \begin{align}t_{\hat y_1}t_{\hat y_l}=t_{\hat y_{l-1}}+t_{\hat y_{l+1}}.\end{align}

As before, {\sl we shall use the symbol $\Box$ for any element in the two-sided ideal $H^{<2323}$}. Let $\xi=q^{\frac12}+q^{-\frac12}$. Recall $*=\{s_0, s_2\}$.

By a simple computation, we get
\begin{align}C_{\hat y_1}=C_\tau C_{s_1}C_{s_2s_3s_2s_1s_0s_2s_3s_2s_3s_0}
=C_\tau C_{s_1}(C_{s_2}C_{{}^\star(s_1s_0s_2s_3s_2s_3s_0)}-C_{{}^*(s_1s_0s_2s_3s_2s_3s_0)}).
  \end{align}

Through a direct (although lengthy) computation we get

\begin{align}C_{\hat y_1}C_{\hat y_l}=&h(C_{\hat y_{l+1}}+C_{\hat y_{l-1}})+\sum\limits_{\substack{y\prec s_3s_2s_1\hat y_l\\ y=s_2s_3s_2s_3y'\ne\tau{}^*\hat y_{l-1}}}\mu(y, s_3s_2s_1\hat y_l)hC_\tau C_{s_1}C_y+\triangle+\Box,\end{align}
where $h=h_{\hat y_0,\hat y_l,\hat y_l}=(q^{\frac12})^4+$lower degree terms  (as a Laurent polynomial in   $q^{\frac12}$), $\triangle$ is a $\mathcal A$-linear combination of some $C_w$ whose coefficients have degrees less than 4 (as a Laurent polynomial in   $q^{\frac12}$).

Therefore, \begin{align}\gamma_{\hat y_1, \hat y_l, \hat y_{l-1}}=\gamma_{\hat y_1, \hat y_l, \hat y_{l+1}}=1.\end{align}

\medskip

Next we   show for any $z\in W$,
$\gamma_{\hat y_1, \hat y_l, z}\ne 0$ if and only if $z=\hat y_{l-1}$ or $\hat y_{l+1}$.

First, by 1.4(f) and 1.4(g), $\gamma_{\hat y_1, \hat y_l, z}\ne 0$ implies $z\in\Gamma_{03}\cap\Gamma_{03}^{-1}$.  By the explicit expressions of the elements in $\Gamma_{03}\cap\Gamma_{03}^{-1}$, we know that $z$ is either  the
first or the third element of a left (resp. right) string with respect to $\{s_2,s_3\}$.

Note that $\hat y_l^{-1}=\hat y_l$ and $(\hat y'_l)^{-1}=\hat y'_l$. Using 1.4(g) and 1.5(d) we get
 \begin{equation}\begin{aligned}
 \gamma_{\hat y_1, \hat y_l, z}&=\gamma_{\hat y_l, z^{-1}, \hat y_1}\ (\text{by\ 1.4(g)}) \\
&=\gamma_{s_2\hat y_l, z^{-1}, s_2\hat y_1}-\gamma_{\hat y'_l, z^{-1}, \hat y_1}\ (\text{by\ 1.5(d)}) \\
&=\gamma_{\hat y_1s_2, s_2\hat y_l, z}-\gamma_{\hat y_1, \hat y'_l, z}\ (\text{by\ 1.4(g)}) \end{aligned}\end{equation}

Since $\gamma_{w, u, v}$ is nonnegative for any $w, u, v\in W$, thus $\gamma_{\hat y_1, \hat y_l, z}\ne 0$ implies  $\gamma_{\hat y_1s_2, s_2\hat y_l, z}\ne 0$.

\begin{itemize}
\item Assume that $z$ is the first element of a left string with respect to $\{s_2,s_3\}$. Then
by 1.5(a) we get $$\gamma_{\hat y_1s_2, s_2\hat y_l, z}=\gamma_{y_1, s_2\hat y_l, s_2z}.$$

When $z$ is the first element of a right string with respect to $\{s_2,s_3\}$, then $z=\hat x_0,\ \hat x'_k,\ \hat y_k$ or $\tilde z_{i,k}^{\star}$ for some $i, k\ge 1$. By 1.5(c'),1.5(f') and Lemma 3.3G we get $\gamma_{y_1, s_2\hat y_l, s_2z}\ne 0$ if and only if  $z=\hat y_{l-1}$ or $\hat y_{l+1}$. And in this case we have $\gamma_{\hat y_1s_2, s_2\hat y_l, \hat y_{l-1}}=\gamma_{\hat y_1s_2, s_2\hat y_l, \hat y_{l+1}}=1$.

When $z$ is the third element of a right string with respect to $\{s_2,s_3\}$, then $z=\hat x_k$ or $\tilde z_{i,k}$ for some $i, k\ge 1$. By 1.5(e') and Lemma 3.3G we get $\gamma_{y_1, s_2\hat y_l, s_2z}=\gamma_{y_1, y_l, s_2zs_3}=0$.

\item Assume that $z$ is the third element of a left string with respect to $\{s_2,s_3\}$. Then
by 1.5(a) we get $$\gamma_{\hat y_1s_2, s_2\hat y_l, z}=\gamma_{y_1, s_2\hat y_l, s_3z}.$$

When $z$ is the first element of a right  string with respect to $\{s_2,s_3\}$, then $z=\tilde x'_k$ or $\tilde z_{i,k}^{\star}$ for some $i, k\ge 1$.  By 1.5(c') and Lemma 3.3G we get $\gamma_{y_1, s_2\hat y_l, s_3z}=\gamma_{y_1, y_l, s_3zs_2}=0$.

When $z$ is the third element of a right string with respect to $\{s_2,s_3\}$, then $z=\tilde x_k, \hat y_k'$ or ${}^\star\tilde z_{i,k}^\star$ for some $i, k\ge 1$.
By 1.5(e'), 1.5(f') and Lemma 3.3G we get $\gamma_{y_1, s_2\hat y_l, s_2z}\ne 0$ if and only if $z=\hat y'_{l-1}$ or $\hat y'_{l+1}$. And in this case $\gamma_{\hat y_1s_2, s_2\hat y_l, \hat y'_{l-1}}=\gamma_{\hat y_1s_2, s_2\hat y_l, \hat y'_{l+1}}=1$.
\end{itemize}

In conclusion, $\gamma_{\hat y_1s_2, s_2\hat y_l, z}\ne 0$ if and only if  $z=\hat y_{l-1}$, $\hat y_{l+1}$,  $\hat y'_{l-1}$ or $\hat y'_{l+1}$. And in each case, $\gamma_{\hat y_1s_2, s_2\hat y_l, z}=1$.

When $z=\hat y'_{l-1}$ or $\hat y'_{l+1}$, i.e. $z=\hat y_{l-1}^\star$ or $\hat y_{l+1}^\star$, by 1.4(k) and (23) we get
\begin{align}\gamma_{\hat y_1, \hat y'_l, z}=\gamma_{\hat y_1, (\hat y'_l)^\star, z^\star}=\gamma_{\hat y_1, \hat y_l, z^\star}=1.\end{align}

Combining formulas (23), (24) and (25), we know $\gamma_{\hat y_1, \hat y_l, z}\ne 0$ if and only if $z=\hat y_{l-1}$ or $\hat y_{l+1}$. Hence (20) is valid.

For $k\ge 2$, we then have    $t_{\hat y_k}t_{\hat y_l}=(t_{\hat y_1}t_{\hat y_{k-1}}-t_{\hat y_{k-2}})t_{\hat y_l}$. Using induction on $k$ we get
$ t_{\hat y_k} t_{\hat y_l} =\sum\limits_{\substack{0\leq i\leq \min \{k, l\}}} t_{\hat y_{k+l-2i}} $.

\medskip

The proof is completed.
 \end{proof}

\medskip

\medskip

{\bf Corollary 3.4G.}
For any $x\in\{\hat y_k, \hat y'_k, \hat z_{i,k}, {}^\star(\hat z_{i,k}), (\hat z_{i,k})^\star, {}^\star(\hat z_{i,k})^\star\ |\  i, k\ge 1\}$ and any
$y\in\{\hat z_{i,k}, {}^\star(\hat z_{i,k}), (\hat z_{i,k})^\star, {}^\star(\hat z_{i,k})^\star\ |\   i, k\ge 1\}$, we have $\pi(t_x)\pi(t_y)=\pi(t_xt_y)$.
\begin{proof}
The result follows from Lemma 3.4D, Corollary 3.4E, Lemma 3.4F and 1.4(k).
\end{proof}

\medskip

 Proposition 3.4A now results from 3.4B-3.4G and 1.4(i).

\medskip

   {\bf Corollary 3.4H.}
 For any left cells $\Gamma,\ \Theta$ in $Y_1$,  we have a bijection
$$\pi: \Gamma \cap\Theta^{-1}\to\{\text{iso. classes of irr. $F$-vector bundles on}\   \{\Theta,\ g\Theta\}\times\{\Gamma,\ g\Gamma\}\}.$$
 \begin{proof}
It follows from Lemma 2.6, subsection 2.4 and Proposition 3.4A.
\end{proof}

\bigskip

  {\bf 3.5.} In this subsection we consider $\Gamma_{013}\cap\Gamma_{013}^{-1}$ and its based ring $J_{\Gamma_{013}\cap\Gamma_{013}^{-1}}.$

  \medskip

  According to  [D, Theorem 6.4], the set $\Gamma_{013}\cap\Gamma^{-1}_{013}$ consists of the elements below:
 \begin{alignat*}{3}
 &\check x_i:=(s_0(^{\sharp *}(\hat x_i)))s_1;\quad &\check x'_k:=(s_0(^{\sharp *}(\hat x'_k)))^{*\sharp}s_0;\quad
&\breve x_i:=s_1\tilde x_is_1;\\
&\breve x'_k:=s_1((({\tilde x'_k})^{*\sharp})s_0);\quad&\check y_k:=(s_0(^{\sharp *}(\hat y_k)))^{*\sharp}s_0;\quad
&\check  y'_k:=s_1((({\hat y'_k})^{*\sharp})s_0);\\
&\check z_{j, k}:=(s_0(^{\sharp *}(\hat z_{j, k})))s_1;\quad
&{}^\star\check z_{j, k},\quad  \check z_{j, k}^\star,\quad & {}^\star\check z_{j, k}^\star.
\end{alignat*}
where   $*=\{s_0, s_2\}$, $\sharp=\{s_1, s_2\}$ and $\star=\{s_2, s_3\},\ i\ge 0,  j, k\geq 1$.

\medskip

  {\bf Proposition 3.5A.}
The bijection
$$\pi: \Gamma_{013}\cap \Gamma^{-1}_{013}\longrightarrow \{\text{iso. classes of irr. $F$-vector bundles on}\ Z_{013}\times Z_{013}\}\subset M_{2\times 2}(R_{F^\circ}),$$

\medskip

$\check x_0 \longmapsto  \left(\begin{array}{cc}
 {\mathbf 1} &0\\
0 &{\mathbf 1}
\end{array}\right),\quad \check x_k \longmapsto
\left(\begin{array}{cc}
\eta^{-k} & 0 \\
0 & \eta^{k}
\end{array}\right),\quad
\check x'_k \longmapsto  \left(\begin{array}{cc}
0 & \eta^{-k} \\
\eta^{k} & 0
\end{array}\right),$

\medskip

$
\breve x_0 \longmapsto  \left(\begin{array}{cc}
0 & {\mathbf 1} \\
{\mathbf 1} & 0
\end{array}\right),\quad
\breve x_k \longmapsto
\left(\begin{array}{cc}
0 & \eta^{k} \\
\eta^{-k} & 0
\end{array}\right),\quad
\breve x'_k \longmapsto  \left(\begin{array}{cc}
\eta^{k} & 0 \\
0 & \eta^{-k}
\end{array}\right),$

\medskip

$
\check y_k \longmapsto  \left(\begin{array}{cc}
V(k) & 0 \\
0 & V(k)
\end{array}\right),\quad
\check y'_k \longmapsto  \left(\begin{array}{cc}
0 & V(k) \\
V(k) & 0
\end{array}\right),$

\medskip

$
\check z_{j,k} \longmapsto   \left(\begin{array}{cc}
\eta^{-j}\otimes V(k) & 0 \\
0 & \eta^j\otimes V(k)
\end{array}\right),\quad
{}^\star\check z_{j,k} \longmapsto   \left(\begin{array}{cc}
0 & \eta^j\otimes V(k)\\
\eta^{-j}\otimes V(k) & 0
\end{array}\right),$

\medskip

$
\check z_{j,k}^\star \longmapsto   \left(\begin{array}{cc}
0& \eta^{-j}\otimes V(k)  \\
\eta^j\otimes V(k) & 0
\end{array}\right),\quad
{}^\star\check z_{j,k}^\star \longmapsto   \left(\begin{array}{cc}
\eta^j\otimes V(k) & 0\\
0& \eta^{-j}\otimes V(k)
\end{array}\right),$

\medskip

\noindent induces a based ring isomorphism
$$\pi: J_{\Gamma_{013}\cap \Gamma^{-1}_{013}}\to  K_{F}(Z_{013}\times Z_{013}).$$

\medskip

We prove the proposition through a series of lemmas.

\medskip

First we establish a bijection between $\Gamma^{-1}_{03}\cap \Gamma_{03}$ and $\Gamma^{-1}_{013}\cap \Gamma_{013}$. Define $\phi: \Gamma^{-1}_{03}\longrightarrow\Gamma^{-1}_{013}$ as follows:
$$\phi(w)=\begin{cases}{s_1w,}&{ \text {if $a(s_1w)=4$,}}\\
{s_0(^{\sharp*}w),}&{\text{otherwise.}}
\end{cases}$$
Then define $\phi': \Gamma_{03}\longrightarrow \Gamma_{013}$ by
$\phi'(w)=(\phi(w^{-1}))^{-1}$. 
Note that for $w\in \Gamma_{03}\cap \Gamma^{-1}_{03}$, we have $\phi(w)\in \Gamma_{03}\cap\Gamma_{013}^{-1}$ and $\phi'(w)\in \Gamma_{013}\cap\Gamma_{03}^{-1}$. We can then define  $\varphi: \Gamma_{03}\cap \Gamma^{-1}_{03}\longrightarrow\Gamma_{013}\cap \Gamma^{-1}_{013}$ by $\varphi(w):=\phi\circ\phi'(w)=\phi'\circ\phi(w)$.

\medskip

{\bf Lemma 3.5B.}
The map $\phi$ is   bijective. In particular,  $\varphi$ is bijective.
\begin{proof}
First, we need to verify $\phi$ is well-defined.
Recall $*=\{s_0,s_2\}, \sharp=\{s_1,s_2\}$ and $\star=\{s_2,s_3\}$.
Let $w=s_0s_3w'\in\Gamma_{03}^{-1}$ for some $w'\in W$ with  $l(w)=l(w')+2$.

If $a(s_1w)=4$, then $s_1w$ is clearly in $\Gamma_{013}^{-1}$.
If $a(s_1w)>4$, we need to show $a(s_0({}^{\sharp*}w))=4$ and $ {L} (s_0({}^{\sharp*}w))=\{s_0, s_1, s_3\}$. Since $w\in\Gamma_{03}^{-1}$ and ${}^*w\in\Gamma_{23}^{-1}$, we know either ${}^*w=s_0w$ or $s_2w$.

\medskip

Note $ s_2\in {L} (s_2w)$ and $s_1\not\in L(s_2w)$. We will show

(a) if $ {L} (s_2w)=\{s_0,s_2\}$, then $s_2w\in\Gamma_{02}^{-1}$ and either $a(s_1w)=4$ or $s_0({}^{\sharp*}w)\in\Gamma_{013}^{-1}$ and $a(s_1w)=9$, in particular ${}^*w=s_0w$;

(b) if $ {L} (s_2w)\ne\{s_0,s_2\}$, then either $s_2w={}^*w\in\Gamma_{23}^{-1}$ and $a(s_1w)=4$ or $ {L} (s_2w)=\{s_0,s_2,s_3\}$ and $s_1w=s_0({}^{\sharp *}w)\in\Gamma_{013}^{-1}$.

\medskip

So we are concerned about  which right cell $s_2w$ is in.

\begin{itemize}
\item Assume $ {L} (s_2w)=\{s_0,s_2\}$. Then $w=s_0s_3s_2s_3s_2\tilde w$ and $s_2w=s_2s_0s_3s_2s_3s_2\tilde w$  with $ {L} (w)=\{s_0, s_3\}$  for some $\tilde w\in W$, and $l(w)=l(\tilde w)+5$, so $s_2w$ and $w$ are in the same left string with respect to $\{s_2,s_3\}$, indicating $a(s_2w)=a(w)=4$. Moreover $s_2w$ is either in $\Gamma_{02}^{-1}$ or $(\Gamma'_{02})^{-1}$.

We claim $s_2w$ is in $\Gamma_{02}^{-1}$.
Otherwise assume $s_2w\in (\Gamma'_{02})^{-1}$, then ${}^\star(s_2w)\in (\Gamma''_2)^{-1}$.
\begin{itemize}
\item If $w$ is the first element of the left string with respect to  $\{s_2,s_3\}$, then $w=s_0s_3s_2s_3s_2\tilde w$ with $s_2\not\in {L} (s_0s_2s_3s_2\tilde w)$ and $s_2w$ is the second element. So ${}^\star(s_2w)=s_2w$ with $s_0\in {L} ({}^\star(s_2w))$, contradicting  the assumption.
\item If $w$ is the second element of the left string with respect to $\{s_2,s_3\}$, then $w=s_0s_3s_2s_3s_2s_0s_2\tilde w'$ for some $\tilde w'\in W$ with $l(w)=l(\tilde w')+7$ and $s_2w$ is the third element. So ${}^\star(s_2w)=s_0s_2s_3s_2s_0s_2\tilde w'$ with $s_0\in {L} ({}^\star(s_2w))$, contradicting  the assumption.
\item If $w$ is the third element of the left string with respect to $\{s_2,s_3\}$, then $w=s_0s_3s_2s_3s_2s_0s_2s_3\dot w$ for some $\dot w\in W$ with $l(w)=l(\dot w)+8$ and $s_2w=w_{023}\dot w$, indicating $a(s_2w)=a(w_{023})=9$, contradicting  the assumption.
\end{itemize}

Therefore, $s_2w$ has to be in $\Gamma_{02}^{-1}$. Moreover, by [D], ${}^\sharp(s_2w)=s_1s_2w\in(\Gamma_{01})^{-1}$ and  ${}^*(s_1s_2w)\in(\Gamma'_{12})^{-1}$.

Next we discuss all the possible reduced expressions of ${}^*(s_1s_2w)$.
\begin{itemize}
\item If ${}^*(s_1s_2w)=s_2s_1s_2w=s_2s_1s_2s_0s_3s_2s_3s_2\tilde w\in (\Gamma'_{12})^{-1}$, then ${}^*(s_1s_2w)=s_1s_2s_1w$. So $
w\underset{LR}\le s_1w\underset{LR}\le {}^*(s_1s_2w)$, by 1.4(b) we get $a(s_1w)=4$.
\item If ${}^*(s_1s_2w)=s_0s_1s_2w=s_1s_2s_0s_3s_2s_3\tilde w\in (\Gamma'_{12})^{-1}$, then $s_1\in {L} (s_0s_3s_2s_3\tilde w)$. Thus $s_0s_3s_2s_3\tilde w=s_0s_3s_2s_3s_1s_2s_3s_2\hat w$ and $s_1w=s_1s_0s_2s_3s_2s_3s_1s_2s_3s_2\hat w=s_0w_{123}\hat w$ for some $\hat w\in W$ with $l(s_1s_2s_3s_2\hat w)=l(\hat w)+4$.
In this case, $a(s_1w)=9$, $s_0({}^{\sharp*}w)=s_0s_3s_2s_1s_3s_2s_3s_2\hat w$ with $\{s_0,s_1,s_3\}\subset L(s_0({}^{\sharp*}w))$ and $s_2w=s_0s_2s_0({}^{\sharp*}w)$. So $s_3s_2s_3s_2\underset{LR}\le s_0({}^{\sharp*}w)\underset{LR}\le s_2w$, by 1.4(b) we get $a(s_0({}^{\sharp*}w))=4$ and $s_0({}^{\sharp*}w)\in\Gamma_{013}^{-1}$.
\end{itemize}

\medskip

In conclusion, when $ {L} (s_2w)=\{s_0,s_2\}$, then $s_2w\in\Gamma_{02}^{-1}$ and either $a(s_1w)=4$ or $s_0({}^{\sharp*}w)\in\Gamma_{013}^{-1}$ and $a(s_1w)=9$. So (a) is valid.

\medskip

\item  Assume $L(s_2w)\not=\{s_0,s_2\}$. Since $s_1\not\in L(s_2w)$, we have $L(s_2w)=\{s_2,s_3\}$ or $L(s_2w)=\{s_0,s_2, s_3\}$.
\begin{itemize}
\item If $L(s_2w)=\{s_2,s_3\}$, then $s_2w={}^*w$, indicating $s_2w\in\Gamma_{23}^{-1}$. So $w=s_0s_3s_2s_0s_3s_2s_3\tilde w$ for some $\tilde w\in W$ with $l(w)=l(\tilde w)+7$.
And ${}^{\sharp*}w=s_1s_2w=s_3s_1s_2s_3s_0s_2s_0s_3s_2\tilde w$ is the first element of the left string with respect to $\{s_2, s_3\}$. Otherwise $s_1\in {L} (s_0s_2s_0s_3s_2\tilde w)$ forces $a(w)=a({}^{\sharp*}w)\ge a(w_{012})=6$, contradicting the choice of $w$.
Hence ${}^\star({}^{\sharp*}w)=s_3s_2s_1s_2w=s_3s_1s_2s_1w$ forces $a(w)\le a(s_1w)\le a({}^{\star\sharp*}w)=a(w)$, indicating $a(s_1w)=4$.

\item If $ {L} (s_2w)=\{s_0,s_2,s_3\}$, then $w=s_0s_3s_2s_3s_2s_0s_2s_3\tilde w$ for some $\tilde w\in W$ with $l(w)=l(\tilde w)+8$
and $s_0({}^{\sharp *}w)=s_1w$ with $ {L} (s_0({}^{\sharp *}w))=\{s_0,s_1,s_3\}$.
Since ${}^{\sharp *}w=s_1s_3s_2s_3s_2s_0s_2s_3\tilde w$ and $s_2({}^{\sharp *}w)=s_2s_1s_3s_2s_3s_2s_0s_2s_3\tilde w$, we know ${}^{\sharp *}(s_2({}^{\sharp *}w))={}^\sharp(s_0s_2({}^{\sharp *}w))$.
If ${}^{\sharp *}(s_2({}^{\sharp *}w))=s_1s_0s_2({}^{\sharp *}w)=s_0s_2s_3s_1s_2s_3s_0s_2s_3\tilde w$ with $s_2\in L({}^{\sharp *}(s_2({}^{\sharp *}w)))$, then $s_2$ has to be in ${L} (\tilde w)$, indicating $a(w)=a(w_{023})=9$, contradicting the choice of $w$.
If ${}^{\sharp *}(s_2({}^{\sharp *}w))=s_2s_0s_2({}^{\sharp *}w)=s_0s_2s_0({}^{\sharp *}w)$, then $a(s_0({}^{\sharp *}w))\le a({}^{\sharp *}(s_2({}^{\sharp *}w)))=a(s_2({}^{\sharp *}w))$.
Next we show $a(s_2({}^{\sharp *}w))=4$, then $a(w)=a({}^{\sharp *}w)\le a(s_0({}^{\sharp *}w))\le a(s_2({}^{\sharp *}w))=4$ forces $a(s_0({}^{\sharp *}w))=4$ and $s_0({}^{\sharp *}w)\in\Gamma_{013}^{-1}$.

 If $s_2({}^{\sharp *}w)=s_2s_1s_3s_2s_3s_2s_0s_2s_3\tilde w$ is not the second element of the left string with respect to $\{s_2, s_3\}$, then $s_1\in {L} (s_2s_0s_2s_3\tilde w)$ forces $a(w)\ge a(s_2s_0s_2s_3\tilde w)\ge a(w_{012})=6$, contradicting the choice of $w$.

Hence $s_3s_2({}^{\sharp *}w)={}^{\star\sharp *}w$ forces $a(w)=a({}^{\sharp *}w)\le a(s_2({}^{\sharp *}w))\le a(s_3s_2({}^{\sharp *}w))=a(w)$. Thus $a(s_2({}^{\sharp *}w))=a(w)=4$ and $s_0({}^{\sharp *}w)\in\Gamma_{013}^{-1}$ consequently.
\end{itemize}
\end{itemize}

\medskip

In conclusion, when $L(s_2w)\not=\{s_0,s_2\}$, then either $s_2w={}^*w\in\Gamma_{23}^{-1}$ and $a(s_1w)=4$ or $ {L} (s_2w)=\{s_0,s_2,s_3\}$ and $s_1w=s_0({}^{\sharp *}w)\in\Gamma_{013}^{-1}$. So (b) is valid.

Therefore, $\phi$ is well-defined.

\medskip

Secondly, we will  prove $\phi$ is  injective.

\begin{itemize}
\item If $\phi(w)=s_1w=\phi(u)=s_1u$, clearly we have $w=u$.
\item  If $\phi(w)=s_0({}^{\sharp *}w)=\phi(u)=s_0({}^{\sharp *}u)$, then ${}^{\sharp*}w={}^{\sharp*}u$. Since left star operations are invertible, we get $w=u$.
\item Assume $\phi(w)=s_1w=\phi(u)=s_0({}^{\sharp *}u)\ne s_1u$ and $s_1w\ne s_0({}^{\sharp *}w)$.

By (a) and (b) we get $u=s_0s_3s_2s_3s_2s_1s_2s_3s_2u'$ for some $u'\in W$ with $l(u)=l(u')+9$ and $\phi(u)=s_0s_3s_2s_3s_1s_2s_3s_2u'=\phi(w)=s_1w$ indicates $w=s_0s_3s_2s_1s_3s_2s_3u'$. Since we assume $s_1w\ne s_0({}^{\sharp *}w)$, by the previous discussion, we get either $s_2w\in\Gamma_{02}^{-1}$ or $s_2w={}^*w\in\Gamma_{23}$.
And the latter case would not happen for the reduced expression of $w$, so $s_2w\in\Gamma_{02}^{-1}$ forces $s_1\in L(u')$ and $a(u)\ge a(w_{123})$, contradicting the choice of $u$.
Thus the assumption  would not occur.
\end{itemize}

In conclusion, $\phi$ is  injective.

\medskip

Next, construct an inverse map of $\phi$, which implies that $\phi$ is surjective. Keep the notations in [D] and subsection 3.1.

Define $\psi: \Gamma^{-1}_{013}\longrightarrow\Gamma^{-1}_{03}$ as follows:
$$\psi(z)=\begin{cases}{s_1z,}&{ \text {if $a(s_1z)=4$ }}\\
{^{*\sharp}(s_0z),}&{\text{otherwise.}}
\end{cases}$$

Note $a(s_1z)\leq a(z)=4$ and $ {L} (s_1z)=\{s_0, s_3\}$.

\begin{itemize}
\item [1.] In fact, $a(s_1z)=4$ if and only if $s_1z\in\Gamma^{-1}_{03}$.
\item[2.] Assume $a(s_1z)<4$.
Recall the notations in [D, subsection 5.2], where { $D_{ij...}, D'_{ij...}$ or  $\hat D'_{ij...}$ denotes the left cell} indexed by the $R$-set $\{s_i, s_j, \cdots\}$ in the two-sided cell $D$ whose $a$-function value is 3;
$\hat E_{ij...}$ denotes the left cell indexed by the $R$-set $\{s_i, s_j, \cdots\}$ in the two-sided cell $\hat E$ whose $a$-function value is 2.
\begin{itemize}
\item[(i) ]
If $a(s_1z)=3$,
and $s_1z\in D^{-1}_{03}\cap\Gamma$,
where $\Gamma\in\{D_0, D_1, D'_2, \hat D'_2, D_{03}, D_{13},$

$ D_{02}, D_{12}\}$. By the reduced expressions of elements in the two-sided cell $D$, we have
$s_1z=s_3s_0(s_2s_3s_1s_0)^ks_2s_0y$ or $\tau s_3s_1(s_2s_3s_1s_0)^{k+1}s_2s_0y$ or $\tau s_3s_1s_2s_1s_0s_2y$ for some $k\geq 0, y\in W$ and $y$ is a product of $\tau$ and simple reflections associated with some right star operations.

For the first two cases, $a(z)=3$,  and for the last case, $a(z)=6$, contradicting the choice of $z$.
\item[(ii) ] If $a(s_1z)=3$ and $s_1z\in D^{-1}_{03}\cap\Gamma'$, where $\Gamma'\in\{ D_{013}, D_{2}, D_{3}\}$.
Then $s_1z=s_3s_0(s_2s_3s_1s_0)^ky$  for some $k\geq 1, y\in W$ and $y$ is a product of $\tau$ and simple reflections associated with some right star operations. In this case $a(z)=3$, contradicting the choice of $z$.

\item[(iii) ] If $a(s_1z)=2$, then $s_1z\in \hat E^{-1}_{03}$.

Hence by the reduced expressions of elements in the two-sided cell $\hat E$, we have $s_1z=s_0s_3(s_2s_1s_3\tau)^ky$ for some $k\geq 0, y\in W$ and $y\in\{\tau, s_2, \tau s_2, s_2s_3, \tau s_2s_3, s_2s_1, \tau s_2s_1\}$.
If $0\leq k\leq 1$, then $a(z)=3$, contradicting the choice of $z$.
If $k\geq 2$,
then $z=s_0s_3s_2s_1s_2s_3s_2s_3s_0(s_2s_1s_3\tau)^{k-2}y$.

In this case, $\psi(z)={}^{*\sharp}(s_0z)=s_0s_3s_2s_3s_2s_1s_2s_3s_2s_0(s_2s_1s_3\tau)^{k-2}y\in \Gamma^{-1}_{03}$.
\end{itemize}
\end{itemize}

In conclusion, if $a(s_1z)=4$, then $\psi(z)=s_1z\in\Gamma^{-1}_{03}$;
if $a(s_1z)<4$, then $a(s_1z)=2$ and $\psi(z)={}^{*\sharp}(s_0z)=s_0s_2s_0z\in \Gamma^{-1}_{03}$.

We have proved that $\psi$ is well-defined and clearly $\psi\circ\phi=1$, so $\psi$ is surjective. In order to prove $\phi$ is bijective, it suffices to show $\psi$ is injective.
\begin{itemize}
\item If $\psi(z)=s_1z=\psi(z')=s_1z'$, then clearly we have $z=z'$.
\item  If $\psi(z)={}^{*\sharp }(s_0z)=\psi(z')={}^{*\sharp }(s_0z')$, then $z=z'$ because left star operations are invertible.
\item Assume $\psi(z)=s_1z=\psi(z')={}^{*\sharp }(s_0z')$.
Then  by 2(iii) in the proof, we have $z'=s_0s_3s_2s_1s_2s_3s_2s_3s_0(s_2s_1s_3\tau)^{k}y$ for $k\geq 2, y\in\{\tau, s_2, \tau s_2, s_2s_3, \tau s_2s_3, s_2s_1, \tau s_2s_1\}$
and $\psi(z')=s_0s_2s_0z'=s_1z$. So $z=s_1s_0s_2s_3s_2s_1s_2s_3s_2s_3s_0(s_2s_1s_3\tau)^{k}y=s_0w_{123}s_0(s_2s_1s_3\tau)^{k}y$, indicating $a(z)=a(w_{123})=9$. This contradicts the choice of $z$. Hence this case would not occur.
\end{itemize}

In conclusion, $\psi$ is injective.

\medskip

We have seen that $\psi$ is bijective, hence its inverse $\phi$ is bijective, and $\phi'$ is  bijective as well. As a consequence, $\varphi$ is bijective.

The proof is completed.
\end{proof}

 \medskip

 {\bf Lemma 3.5C.}
If $w\in\Gamma^{-1}_{03}$, then
$C_{s_1}C_w=C_{\phi(w)}+\Box$, where $\Box\in H^{<2323}$.

\begin{proof}
Let $w=s_0s_3w'$ for some $w'\in W$ with $l(w)=l(w')+2$. Since $w\in\Gamma_{03}^{-1}$, we have $s_1w'>w'$.

As before, {\sl we  use the symbol $\Box$ for any element in the two-sided ideal $H^{<2323}$}, and  $\xi=q^{\frac12}+q^{-\frac12}$.

By formula in 1.1(a) and a simple computation we have
$$C_{s_1}C_w=C_{s_1w}+\sum\limits_{\substack{y\prec w\\ L(y)=\{s_0,s_1,s_3\}}}\mu(y, w)C_y.$$

Assume $y\prec w$, $s_0y<y, s_3y<y$ and $s_1y<y$. Let $y=s_0s_1s_3y'$ for some $y'\in W$ with $l(y)=l(y')+3$.

Recall $*=\{s_0,s_2\}$. By 1.3(b), $\mu(y, w)=\tilde\mu({}^*y, {}^*w).$ There are two cases for ${}^*y$.

(a)  Assume ${}^*y=s_1s_3y'$. Then $s_2\in {L} ({}^*y)$ indicates ${}^*y=w_{123}\tilde y$ and $y=s_0w_{123}\tilde y$ for some $\tilde y\in W$ with $l(y)=l(\tilde y)+10$. Moreover, $a(y)=a(w_{123})=9$ and $C_y\in H^{<2323}$.

(b)  Assume ${}^*y=s_2y$. Then
by [D],
we know ${}^*w\in\Gamma_{23}^{-1}$ and there are two cases for ${}^*w$.

\begin{itemize}
\item Assume ${}^*w=s_0w\in\Gamma_{23}^{-1}$. Then $w=s_0s_3s_2s_3s_2\tilde w$ for some $w\in W$ and $w'=s_2s_3s_2\tilde w$ for some $\tilde w\in W$ with $l(w')=l(\tilde w)+3$.
Then $\mu(y, w)=\mu({}^*y, {}^*w)=\tilde\mu(s_2y, s_0w).$

By a  simple discussion we get that $^{*}w=s_0w$ together with ${}^*y=s_2y$ induces either $C_y\in H^{<2323}, a(s_1w)=4$ or $w=s_0s_2s_3s_2s_3s_1s_2s_3s_2\hat w$ for $\hat w\in W$ satisfying $l(w)=l(\hat w)+9$ and $y=\phi(w)=s_0({}^{\sharp *}w)$ with $C_{s_1w}\in H^{<2323}$.

\item Assume ${}^*w=s_2w\in\Gamma_{23}^{-1}$. Then $w=s_0s_3s_2s_0s_3s_2s_3\tilde w$ for some $\tilde w\in W$ with $s_0\tilde w>\tilde w$ and $w'=s_2s_0s_3s_2s_3\tilde w$ for some $\tilde w\in W$ with $l(w')=l(\tilde w)+5$.
Then $\mu(y, w)=\mu({}^*y, {}^*w)=\mu(s_2y, s_2w)$.

By a simple analysis we get that $^{*}w=s_2w$ together with ${}^*y=s_2y$ induces $C_y\in H^{<2323}$. And in this case $a(s_1w)=4$.
\end{itemize}

\medskip

Therefore, $$C_{s_1}C_w=C_{s_1w}+\delta C_{s_0({}^{\sharp *}w)}+\Box=C_{\phi(w)}+\Box,$$ where $\Box\in H^{<2323}$ and $\delta=\begin{cases}{0,}&{\text{if}\ \phi(w)=s_1w}\\{1,}&{\text{if}\ \phi(w)\not=s_1w}\end{cases}$.
Moreover, by the proof of Lemma 3.5B, when $\delta=1$, we have $a(s_1w)=9$.

\medskip

The proof is completed.
\end{proof}

\medskip

{\bf Lemma 3.5D.}
For any $w, u, z\in  \Gamma^{-1}_{03}, v\in W$, we have
\begin{itemize}
\item [(a) ] $\gamma_{\phi(w), v, \phi(z)}=\gamma_{w, v, z}$;
\item [(b) ] $\gamma_{\phi'(w^{-1}), \phi(u), z}=\gamma_{w^{-1}, u, z}$;
\item [(c) ] { $\gamma_{w, \phi'(u^{-1}), \phi'(z^{-1})}=\gamma_{w, u^{-1}, z^{-1}}$.}
\end{itemize}
\begin{proof}
As before,  let $\xi=q^{\frac12}+q^{-\frac12}$.

By Lemma 3.5C, for any $w, z\in  \Gamma^{-1}_{03}, v\in W$, we have $C_{s_1}C_wC_v\in C_{\phi(w)}C_v+H^{<2323}$. And $C_wC_v\in\sum \limits_{\substack{z\in\Gamma^{-1}_{03},\\ s_1z>z\\ \gamma_{w, v, z}\ne0}}hC_z+H^{<2323},$ where $h\in\mathbb Z[\xi]$ is a polynomial in $\xi$ with degree $\le$ 4.
Therefore,  $\gamma_{\phi(w), v, \phi(z)}=\gamma_{w, v, z}$ and part (a) is proved.

Part (b) follows from part (a) and 1.4(g),
since $\gamma_{\phi'(w^{-1}), \phi(u), z}=\gamma_{\phi(u), z^{-1}, (\phi'(w^{-1}))^{-1}} =\gamma_{\phi(u), z^{-1}, \phi(w)} =\gamma_{u, z^{-1}, w}=\gamma_{w^{-1}, u, z}$.

{ Note $\phi'(x)=(\phi(x^{-1}))^{-1}$, for any $x\in\Gamma^{-1}_{03}$. Then part (c) follows from part (b) and 1.4(g):

$\gamma_{w, \phi'(u^{-1}), \phi'(z^{-1})}=\gamma_{w, \phi'(u^{-1}), (\phi(z))^{-1}}=\gamma_{\phi'(u^{-1}), \phi(z), w^{-1}}=\gamma_{u^{-1}, z, w^{-1}}=\gamma_{w,u^{-1},z^{-1}}$.}

The proof is completed.
\end{proof}

The following lemma tells us that $\varphi$ preserves the $\gamma$-relations.

\medskip

{\bf Corollary 3.5E.}
For any $w, z, u\in \Gamma_{03}\cap \Gamma^{-1}_{03}$,
$\gamma_{\varphi(w), \varphi(z), \varphi(u)}=\gamma_{w, z, u}$. In particular, the based ring $J_{\Gamma_{013}\cap \Gamma^{-1}_{013}}$ is isomorphic to $J_{\Gamma_{03}\cap \Gamma^{-1}_{03}}$, { the based module $J_{\Gamma\cap \Gamma^{-1}_{013}}$ is isomorphic to $J_{\Gamma\cap \Gamma^{-1}_{03}}$, and $J_{\Gamma_{013}\cap\Gamma^{-1}}$ is isomorphic to $J_{ \Gamma_{03}\cap\Gamma^{-1}}, $ for any left cell $\Gamma$.}

\begin{proof}
It is deduced from  Lemma 3.5D.
\end{proof}

Combining Corollary 3.5E and Proposition 3.4A, we see that Proposition 3.5A is true.

\medskip

   {\bf Corollary 3.5F.}
   For any left cells $\Gamma,\ \Theta$ in $Y_2$, we have a bijection
 $$\pi: \Gamma \cap\Theta^{-1}\to\{\text{iso. classes of irr. $F$-vector bundles on}\ \{\Theta,\ g\Theta\}\times\{\Gamma,\ g\Gamma\}\}.$$

 \begin{proof}
It follows from Lemma 2.6, subsection 2.4 and Proposition 3.5A.
\end{proof}

\bigskip

 {\bf 3.6.} In this subsection we establish a bijection
  $$\pi: \Gamma_{02}\cap\Gamma_{03}^{-1}\to\{\text{iso. classes of irr. $F$-vector bundles on}\  Z_{03}\times \{\Gamma_{02}\}\}.$$
  It also gives rise to a bijection
  $$\pi': \Gamma_{03}\cap\Gamma_{02}^{-1}\to\{\text{iso. classes of irr. $F$-vector bundles on}\ \{\Gamma_{02}\}\times Z_{03}\}.$$

{
Recall that $\phi': \Gamma_{03}\longrightarrow \Gamma_{013}$ is a bijection.
So is $(\phi')^{-1}: \Gamma_{013}\to\Gamma_{03}$. Then
 $\Gamma_{013}\cap\Gamma_{02}^{-1}\overset{(\phi')^{-1}}{\longrightarrow}\Gamma_{03}\cap\Gamma_{02}^{-1}, \Gamma_{013}\cap\Gamma_{03}^{-1}\overset{(\phi')^{-1}}{\longrightarrow}\Gamma_{03}\cap\Gamma_{03}^{-1} $ are also bijections.

 By Corollary 3.5E, $\pi'$ gives rise to a bijection
  $$\hat\pi': \Gamma_{013}\cap\Gamma_{02}^{-1}\to\{\text{iso. classes of irr. $F$-vector bundles on}\ \{\Gamma_{02}\}\times Z_{013}\}: x\mapsto \pi'((\phi')^{-1}(x)).$$}

Following [D, Theorem 6.4], we have
 $$\Gamma_{02}\cap\Gamma^{-1}_{03}=\{s_2x_k, s_3x_j, s_2y_k, s_2z_{i,k}, s_3z_{i,k}\ |\ i, k\geq 1, j\ge 0\}.$$ Let $*=\{s_0, s_2\}$, $\sharp=\{s_1, s_2\}$ and $\star=\{s_2, s_3\},  i, k\geq 1, j\ge 0$.

For expressing the proposition below, we denote the map in Proposition 3.3A, Proposition 3.4A, and Proposition 3.5A by $\pi_1, \pi_2, \pi_3$ respectively. Then by Corollary 3.5E we get a bijection $\hat\pi_2: \Gamma_{013}\cap\Gamma^{-1}_{03}\mapsto \{\text{iso.classes of irr. $F$-vector bundles on}\ Z_{03}\times Z_{013}\}: x\mapsto \pi_2((\phi')^{-1}(x)).$

  \medskip

  {\bf Proposition 3.6A.}
There is a bijection
$$\pi: \Gamma_{02}\cap \Gamma^{-1}_{03}\to\{\text{iso. classes of irr. $F$-vector bundles on}\  Z_{03}\times\{\Gamma_{02}\}\}\subset M_{2\times1}(R_{F^\circ}),$$

\medskip

$s_2x_k\longmapsto \left(\begin{array}{c}
\eta^{-k}  \\
\eta^{k}
\end{array}\right),
s_3x_k\longmapsto  \left(\begin{array}{c}
\eta^{k}  \\
\eta^{-k}
\end{array}\right),$

\medskip

$s_2y_k\longmapsto   \left(\begin{array}{c}
V(k)  \\
V(k)
\end{array}\right),$

\medskip

$s_2z_{i,k}\longmapsto  \left(\begin{array}{c}
\eta^{-i}\otimes V(k)\\
\eta^{i}\otimes V(k)\end{array}\right),
s_3z_{i,k}\longmapsto  \left(\begin{array}{c}
\eta^{i}\otimes V(k)\\
\eta^{-i}\otimes V(k)\end{array}\right)$

\medskip

\noindent induces a bijection
$$\pi: J_{\Gamma_{02}\cap \Gamma^{-1}_{03}}\to  \{\text{iso. classes of irr. $F$-vector bundles on}\ Z_{03}\times\{\Gamma_{02}\}\}, t_x\mapsto \pi(x)$$

such that take any $x\in\Gamma_{02}\cap\Gamma_{03}^{-1}$,
we have

(a) $\pi(t_xt_y)=\pi(t_x)\pi_1(t_y)$, for any
$y\in\Gamma_{02}\cap\Gamma_{02}^{-1}$;

(b) $\pi_2(t_xt_y)=\pi(t_x)\pi'(t_y)$, for any $y\in\Gamma_{03}\cap\Gamma_{02}^{-1}$;

(c) $\hat\pi_2(t_xt_y)=\pi(t_x)\hat\pi'(t_y)$,  for any $y\in\Gamma_{013}\cap\Gamma_{02}^{-1}$.

\bigskip

Clearly $\pi$ is bijective.   We prove the proposition through a series of lemmas.

\medskip

{\bf Case 1: } Let  $x\in\Gamma_{02}\cap\Gamma_{03}^{-1}$ and $y$ be in $\Gamma_{02}\cap\Gamma_{02}^{-1}$.

\medskip

{\bf Lemma 3.6B.}
For $x\in\{s_2x_i, s_3x_k\ |\ i\ge 1, k\ge 0\}, y\in\{d_{02}, x_k\ |\ k\ge 0\}$, we have $\pi(t_xt_y)=\pi(t_x)\pi_1(t_y)$.
\begin{proof}
First of all, when $y=d_{02}$, by Lemma 3.3C, we get $y$ is a distinguished involution, then for any $x\in\{s_2x_i, s_3x_k\ |\ i\ge 1, k\ge 0\}$, $t_xt_y=t_x$ and the result is trivial.

Using 1.4(k) and  the fact $x_0=d_{02}^\star, s_2x_k=(s_2x_k)^\star, s_3x_j=(s_3x_j)^\star$, when $y=x_0$,  we get for any $x\in\{s_2x_i, s_3x_k\ |\ i\ge 1, k\ge 0\}$, $t_xt_y=t_{x^\star}$ and the result is true.

Set $t_{s_2x_j}=t_{s_3x_m}=0$ for any integers $j\le 0, m<0$. Next we prove for any $k\ge 0, l\ge 1$, we have \begin{align}t_{s_3x_k}t_{x_l}=t_{s_3x_{k-l}}+t_{s_3x_{k+l}}+t_{s_2x_{l-k}}.\end{align}

Note $s_3x_{k}={}^{\star*\sharp}(\tau s_3x_{k-1})$, $s_2x_l={}^{\star*\sharp}(\tau s_2x_{l+1})$ and $s_3x_0={}^{\star*\sharp}(\tau s_2x_1)$ for any $k, l \ge 1$. Then using 1.4(j) and 1.4(k) we get $\gamma_{s_3x_k, x_l, p}=\gamma_{{}^{\star*\sharp}(\tau s_3x_{k-1}), x_l, p}=\gamma_{s_3x_{k-1}, x_l, \tau({}^{\sharp*\star} p)}$ for any $p\in W$. So by induction it suffices to prove 
for any $l\ge1$ we have , \begin{align} t_{s_3x_0}t_{x_l}=t_{s_3x_{l}}+t_{s_2x_l}.\end{align}

To show (27),  we need to determine which $m\in\Gamma_{02}\cap\Gamma_{03}^{-1}$ satisfies $\gamma_{s_3x_0, x_l, m}\ne 0$. Note $m$ is in some left string with respect to $\{s_2, s_3\}$.

\begin{itemize}
\item Assume $m$ is the first element of the left string with respect to $\{s_2, s_3\}$, then by Lemma 1.5(a) and Lemma 3.3B we get $\gamma_{s_3x_0, x_l, m}=\gamma_{x_0, x_l, s_2m}\ne 0$ if and only if $s_2m=x_l$. So \begin{align}\gamma_{s_3x_0, x_l, s_2x_l}=1.\end{align}

\item Assume $m$ is the second element of the left string with respect to $\{s_2, s_3\}$, then by Lemma 1.5(b) and Lemma 3.3B we get $\gamma_{s_3x_0, x_l, m}=\gamma_{x_0, x_l, s_3m}+\gamma_{x_0, x_l, s_2m}\ne 0$ if and only if $s_3m=x_l$ or $s_2m=x_l$. This contradicts the assumption. So $m$ being the second element of the left string with respect to $\{s_2, s_3\}$ would not occur.

\item Assume $m$ is the third element of the left string with respect to $\{s_2, s_3\}$, then by Lemma 1.5(c) and Lemma 3.3B we get $\gamma_{s_3x_0, x_l, m}=\gamma_{x_0, x_l, s_3m}\ne 0$ if and only if $s_3m=x_l$. So \begin{align}\gamma_{s_3x_0, x_l, s_3x_l}=1.\end{align}
\end{itemize}

Combining formulas (28) and (29), we get formula (27) is valid.
Then we get (26) inductively.

Using 1.4(k) and $s_2x_i={}^\star(s_3x_i)$ for any $i\ge 1$, we can get
$$\pi(t_{s_2x_k})\pi_1(t_{x_l})=\pi(t_{s_2x_{k-l}})+\pi(t_{s_2x_{k+l}})+\pi(t_{s_3x_{l-k}}).$$

\medskip

The proof is completed.
\end{proof}

\medskip

{\bf Lemma 3.6C.}
For $x\in\{s_2x_i, s_3x_k\ |\ i\ge 1, k\ge 0\}, y\in\{y_i, z_{i,j}\ |\ i, j\ge1\}$, we have $\pi(t_xt_y)=\pi(t_x)\pi_1(t_y)$.
\begin{proof}
First of all, note $s_2x_i={}^\star(s_3x_i), s_2z_{i,k}={}^\star(s_3z_{i,k})$ for any $i, k\ge 1$ with $\star=\{s_2,s_3\}$.  Set $s_3z_{0,k}=s_3\dot y_k=s_2y_k$. Then using 1.4(k), it suffices to prove for any $j\ge 0, k\ge 1$, we have

\begin{align}t_{s_3x_j}t_{y_k}=t_{s_3x_j}t_{\dot y_k}=t_{s_3z_{j,k}}.\end{align}

To prove (30), note $\dot y_k=y_k^\star,  s_3z_{j,k}=(s_3z_{j,k})^\star,  s_3x_i={}^{\star*\sharp}(\tau s_3x_{i-1})$ and $s_3z_{i,k}={}^{\star*\sharp}(\tau s_3z_{i-1,k})$ for any $i, k\ge 1, j\ge 0$, thus by induction it  suffices to show

\begin{align}
t_{s_3x_0}t_{y_k}=t_{s_3z_{0,k}}=t_{s_2y_k}.\end{align}

To prove (31), similar to the proof of Lemma 3.6B, using 1.5 and Lemma 3.3C we know for any $p\in W$ such that $\gamma_{s_3x_0, y_k, p}\ne 0$, $p$ has to be $s_2y_k$
and $\gamma_{s_3x_0, y_k, s_2y_k}=1.$
 Thus (31) is valid.

Using 1.4(k) and $s_2x_k={}^\star(s_3x_k)$ for any $k\ge 1$, we can get for any $i, k\ge 1$,
$\pi(t_{s_2x_i})\pi_1(t_{y_k})=\pi(t_{s_2x_i})\pi_1(t_{\dot y_k})=\pi(t_{s_2z_{i,k}}).$

Along with Lemma 3.3E and Lemma 3.6B, we get for any $i, l, k\ge 1$,
$\pi(t_{s_2x_i})\pi_1(t_{z_{l,k}})=\pi(t_{s_2z_{i-l,k}})+\pi(t_{s_2z_{i+l,k}})+\pi(t_{s_3z_{l-i,k}}),$
$\pi(t_{s_3x_i})\pi(t_{z_{l,k}})=\pi(t_{s_3z_{i-l,k}})+\pi(t_{s_3z_{i+l,k}})+\pi(t_{s_2z_{l-i,k}}).$

\medskip

The proof is completed.
\end{proof}

\medskip

{\bf Lemma 3.6D.}
For $x\in\{s_2y_i\ |\ i\ge 1\}, y\in\{d_{02}, x_k, y_i, z_{i,j}\ |\ i, j\ge1, k\ge 0\}$, we have $\pi(t_xt_y)=\pi(t_x)\pi_1(t_y)$.

\begin{proof}
First of all, when $y=d_{02}$ or $x_0$ is trivial.

Let $i, k, l$ be positive intergers. Note $s_2y_0=s_3x_0$. Set $z_{i,0}=x_i$. Then it suffices to prove (a) $\pi(t_{s_2y_k})\pi_1(t_{x_i})=\pi(t_{s_2z_{i,k}})+\pi(t_{s_3z_{i,k}});$ (b) $\pi(t_{s_2y_k})\pi_1(t_{y_l})=\pi(t_{s_2y_k})\pi_1(t_{\dot y_l})=\sum\limits_{\substack{0\leq j\leq \min \{k, l\}}}\pi(t_{s_2y_{k+l-2j}});$

For any $m\in\Gamma_{02}\cap\Gamma_{03}^{-1}$ such that $\gamma_{s_2y_k, x_i, m}\ne 0$, $m$ has to be in the left string with respect to $\{s_2, s_3\}$.
Note $s_2y_k$ is the second element of the left string with respect to $\{s_2, s_3\}$. Similar to the proof of Lemma 3.6B, by 1.5 and Lemma 3.3E, we get
$m$ has to be $s_2z_{i,k}$ or $s_3z_{i,k}$ and
$$\gamma_{s_2y_k, x_i, s_2z_{i,k}}=\gamma_{s_2y_k, x_i, s_3z_{i,k}}=1.$$

Then we get  (a) directly.

\medskip

Note $s_2y_0=s_3x_0$. For (b), we first prove for any $l\ge 1$,
\begin{align} t_{s_2y_1}t_{y_l}=t_{s_2y_{l-1}}+t_{s_2y_{l+1}}.\end{align}

To show (32), we know for any $p\in\Gamma_{02}\cap\Gamma_{03}^{-1}$ such that $\gamma_{s_2y_1, y_l, p}\ne 0$, $p$ has to be in the left string with respect to $\{s_2, s_3\}$. Similarly, by 1.5 and Lemma 3.3G, we get $p$ has to be $s_2y_{l+1}$ or $s_2y_{l-1}$ and
$\gamma_{s_2y_1, y_l, s_2y_{l+1}}=\gamma_{s_2y_1, y_l, s_2y_{l-1}}=1.$

Thus (32) is valid.

\medskip

And $t_{s_2y_k}t_{y_l}=(t_{s_2y_1}t_{y_{k-1}}-t_{s_2y_{k-2}})t_{y_l}$, so by induction on $k$, Lemma 3.3G,  part (a), $s_2y_m=(s_2y_m)^\star$ and 1.4(k), we get part (b).

Along with Lemma 3.3E and 1.4(i), we can get
$$\pi(t_{s_2y_k})\pi_1(t_{z_{i,l}})=\sum\limits_{\substack{0\leq j\leq \min \{k, l\}}}\pi(t_{s_2z_{i,k+l-2j}})+\pi(t_{s_3z_{i,k+l-2j}}).$$

The proof is completed.
\end{proof}

\medskip

{\bf Corollary 3.6E.}
For $x\in\{s_2z_{i,j}, s_3z_{i,j}\ |\ i, j\ge 1\}$, $y\in\Gamma_{02}\cap\Gamma_{02}^{-1}$, we have $\pi(t_xt_y)=\pi(t_x)\pi_1(t_y)$.

\begin{proof}
 Note ${}x_0=(d_{02})^\star,  (s_2z_{i,k})^\star=s_2z_{i,k}, (s_3z_{i,k})^\star=s_3z_{i,k}, s_3z_{i,k}={}^\star(s_2z_{i,k})$.Then the result follows  Lemma 3.6C, Lemma 3.3E, Lemma 3.3G , Corollary 3.3H and 1.4(k).
\end{proof}

\medskip

{ Combining Lemma 3.6B - Lemma 3.6E,  we see that assertion (a) in Proposition 3.6A is true.}

\medskip

{\bf Case 2: } Let $x\in\Gamma_{02}\cap\Gamma_{03}^{-1}$ and $y$ be in $\Gamma_{03}\cap\Gamma_{02}^{-1}$.

{\bf Lemma 3.6F.}
For $x, y^{-1}\in\{s_2x_i, s_3x_k\ |\ i\ge 1, k\ge 0\}$,  we have $\pi_2(t_xt_y)=\pi(t_x)\pi'(t_y)$.
\begin{proof}
Note$(s_2x_i)^{-1}=((s_3x_i)^{-1})^\star, s_2x_i={}^\star(s_3x_i)$, where $\star=\{s_2,s_3\}$.
Let $k, j$ be any non-negative integers. Set $t_{\hat x_m}=t_{\tilde x'_0}=t_{\tilde x'_m}=0$ for any $m<0$.
Then using 1.4(k) it suffices to prove \begin{align}t_{s_3x_k}t_{(s_3x_j)^{-1}}=
t_{\tilde x'_{k-j}}+t_{\tilde x_{k+j}}+t_{\hat x_{j-k}}.\end{align}

Since $s_3x_k=\tau {}^{\star\sharp*}(s_3x_{k-1}), \hat x_{k}=\tau {}^{\star\sharp*}(\hat x_{k+1}), \tilde x_k=\tau {}^{\star\sharp*}(\tilde x_{k-1})$, $\tilde x'_k=\tau {}^{\star\sharp*}(\tilde x'_{k-1})$ and $\tilde x'_1=\tau {}^{\star\sharp*}(\hat x_0)$ for any $k\ge 1$. Then using 1.4(j) and 1.4(k) we only need to prove for any $j\ge 0$,
\begin{align}t_{s_3x_0}t_{(s_3x_j)^{-1}}=
t_{\tilde x_{j}}+t_{\hat x_{j}}.\end{align}

For (34), we consider $p\in\Gamma_{03}\cap\Gamma_{03}^{-1}$ such that $\gamma_{s_3x_0, (s_3x_j)^{-1}, p}\ne 0$.
 According to subsection 3.4 we know that for any $p\in\Gamma_{03}\cap\Gamma_{03}^{-1}$, $p$ has to be the first or third element in a left string with respect to $\{s_2, s_3\}$. And $s_3x_0$ is a second element in a left string with respect to $\{s_2, s_3\}$.
Similar to the proof of Lemma 3.6B,
using 1.5(a), 1.5(c), 1.4(i) and Lemma 3.6B, we get $p=\hat x_j$ or $\tilde x_j$ and  $\gamma_{s_3x_0, (s_3x_j)^{-1}, \hat x_j}=\gamma_{s_3x_0, (s_3x_j)^{-1}, \tilde x_j}=1.$ so (34) is valid.

\medskip

Moreover, we can get

$$\pi(t_{s_2x_k})\pi'(t_{(s_2x_l)^{-1}})=
\pi_2(t_{\hat x_{k-l}})+\pi_2(t_{\hat x'_{k+l}})+\pi_2(t_{\tilde x'_{l-k}}),$$
$$\pi(t_{s_3x_i})\pi'(t_{(s_2x_l)^{-1}})=
\pi_2(t_{\tilde x_{i-l}})+\pi_2(t_{\tilde x'_{i+l}})+
\pi_2(t_{\hat x'_{i-k}});$$
$$\pi(t_{s_2x_k})\pi'(t_{(s_3x_j)^{-1}})=\pi_2(t_{\hat x'_{k-j}})+\pi_2(t_{\hat x_{k+j}})+\pi_2(t_{\tilde x_{j-k}}).$$

\medskip

The proof is completed.
\end{proof}

\medskip

{\bf Lemma 3.6G.}
For $x\in\{s_2x_i, s_3x_k\ |\ i\ge 1, k\ge 0\}, y^{-1}\in\{s_2y_i, s_2z_{i,j}, s_3z_{i,j}|\ i, j\ge 1\}$,  we have $\pi_2(t_xt_y)=\pi(t_x)\pi'(t_y)$.

\begin{proof}
We first prove for any $j\ge 0, k\ge 1$, \begin{align}t_{s_3x_j}t_{(s_2y_k)^{-1}}=t_{{}^\star(\hat z_{j,k})}+t_{{}^\star(\hat z_{j,k})^\star}.\end{align}

Note $s_3x_i=\tau{}^{\star\sharp*}(s_3x_{i-1})$ for any $i\ge 1$. Then using 1.4(j) and 1.4(k), to show (35) it suffices to prove \begin{align}t_{s_3x_0}t_{(s_2y_k)^{-1}}=t_{\hat y'_k}+t_{\hat y_k}.\end{align}

 We know that for any $p\in\Gamma_{03}\cap\Gamma_{03}^{-1}$ such that $\gamma_{s_3x_0, (s_2y_k)^{-1}, p}\ne 0$, $p$ has to be the first or third element in a left string with respect to $\{s_2, s_3\}$, while $s_3x_0$ is a second element of a left string with respect to $\{s_2, s_3\}$.
Similar to the proof of Lemma 3.6B, using 1.5(a), 1.5(c) and Lemma 3.6D we can get
 $\gamma_{s_3x_0, (s_2y_k)^{-1}, p}\ne 0$ if and only if
$p=\hat y'_k$ or $\hat y_k$. And $\gamma_{s_3x_0, (s_2y_k)^{-1}, \hat y'_k}=\gamma_{s_3x_0, (s_2y_k)^{-1}, \hat y'_k}=1$. So (35) is valid.

Note ${}^\star(\hat z_{j, k})=\tau{}^{\star\sharp*}({}^\star(\hat z_{j-1,k})), {}^\star(\hat z_{j, k})^\star=\tau{}^{\star\sharp*}(^\star(\hat z_{j, k})^\star)$
Then we get (35) inductively by 1.4(j) and 1.4(k).

Note $s_2x_i={}^\star(s_3x_i)$ with $\star=\{s_2, s_3\}$ for any $i\ge 1$, then using 1.4(k) we can get $$\pi(t_{s_2x_i})\pi'(t_{(s_2y_k)^{-1}})=\pi_2(t_{\hat z_{i,k}})+\pi_2(t_{(\hat z_{i,k})^\star}).$$

\medskip

Set $t_{\hat z_{m, k}}=t_{{}^\star(\hat z_{m, k})}=t_{(\hat z_{m, k})^\star}=t_{{}^\star(\hat z_{m, k})^\star}=0$ for any $m<0$ and $s_3z_{0,k}=s_2y_k$ for any $k\ge 1$. Then ${}^\star(\hat z_{0,k})=\hat y'_k$ and ${}^\star(\hat z_{0,k})^\star=\hat y_k.$

Similarly, we can show
$$\begin{aligned}&\pi(t_{s_2x_i})\pi'(t_{(s_2z_{l,k})^{-1}})=\pi_2(t_{(\hat z_{i+l,k})^\star})+\pi_2(t_{\hat z_{i-l,k}})+\pi_2(t_{{}^\star(\hat z_{l-i,k})^\star}),\\
&\pi(t_{s_3x_j})\pi'(t_{(s_2z_{l,k})^{-1}})=\pi_2(t_{{}^\star(\hat z_{j+l,k})^\star})+\pi_2(t_{{}^\star(\hat z_{j-l,k})})+\pi_2(t_{(\hat z_{l-j,k})^\star}),\\
&\pi(t_{s_2x_i})\pi'(t_{(s_3z_{l,k})^{-1}})=\pi_2(t_{\hat z_{i+l,k}})+\pi_2(t_{(\hat z_{i-l,k})^\star})+\pi_2(t_{{}^\star(\hat z_{l-i,k})}),\\
&\pi(t_{s_3x_j})\pi'(t_{(s_3z_{l,k})^{-1}})=\pi_2(t_{{}^\star(\hat z_{j+l,k})})+\pi_2(t_{{}^\star(\hat z_{j-l,k})^\star})+\pi_2(t_{\hat z_{l-j,k}}).\end{aligned}$$

\medskip

The proof is completed.
\end{proof}

\medskip

{\bf Lemma 3.6H.}
For $x\in\{s_2y_i\ |\ i\ge 1\}, y^{-1}\in\{s_2y_i, s_2z_{i,j}, s_3z_{i,j}\ |\ i, j\ge 1\}$,  we have $\pi_2(t_xt_y)=\pi(t_x)\pi'(t_y)$.
\begin{proof}
Set $z_{k,0}=x_k$ for any $k\ge 1$ and $\hat y_0=\hat x_0, \hat y'_0=\tilde x_0$. Then $\hat z_{k,0}=\hat x_k$.

We first prove for any $k, l\geq 1$,
\begin{align}\pi(t_{s_2y_k})\pi(t_{(s_2y_l)^{-1}})
=\sum\limits_{\substack{0\leq j\leq \min \{k, l\}}}\pi(t_{\hat y_{k+l-2j}})+\pi(t_{\hat y'_{k+l-2j}}).\end{align}

Since by 1.4(i) and Lemma 3.6D, we get  $t_{s_2y_k}t_{(s_2y_l)^{-1}}=t_{s_2y_k}(t_{y_{l-1}}t_{(s_2y_1)^{-1}}-t_{(s_2y_{l-2})^{-1}})$.

Hence to prove (37), it suffices to show
\begin{align} t_{s_2y_k}t_{(s_2y_1)^{-1}}=t_{\hat y_{k+1}}+t_{\hat y_{k-1}}+t_{\hat y'_{k+1}}+t_{\hat y'_{k-1}}.\end{align}

 To prove (38), we know that for any $p\in\Gamma_{03}\cap\Gamma_{03}^{-1}$ such that $\gamma_{s_2y_k, (s_2y_1)^{-1}, p}\ne 0$, $p$ has to be the first or third element in a right string with respect to $\{s_2, s_3\}$, while $(s_2y_k)^{-1}$ is a second element in a right string with respect to $\{s_2, s_3\}$. Using 1.5(c'),1.5(e') and Lemma 3.6D we get
 $\gamma_{s_2y_k, (s_2y_l)^{-1}, p}\ne 0$ if and only if
$p=\hat y_{k-1}$, $\hat y_{k+1}$, $\hat y'_{k-1}$ or $\hat y'_{k+1}$. And $$\gamma_{s_2y_k, (s_2y_1)^{-1}, \hat y_{k-1}}=\gamma_{s_2y_k, (s_2y_1)^{-1}, \hat y_{k+1}}=\gamma_{s_2y_k, (s_2y_1)^{-1}, \hat y'_{k-1}}=\gamma_{s_2y_k, (s_2y_1)^{-1}, \hat y'_{k+1}}=1.$$

Then we get (38) directly.

\medskip

By Lemma 3.6B, Lemma 3.6D, Lemma 3.6F, Lemma 3.6G and 1.4(i), we get
$$\begin{aligned}&\pi(t_{s_2y_k})\pi(t_{(s_2z_{i,l})^{-1}})=\sum\limits_{\substack{0\leq j\leq \min \{k, l\}}}\pi(t_{{}^\star(\hat z_{i,k+l-2j})^\star})+\pi(t_{(\hat z_{i,k+l-2j})^\star});\\
&\pi(t_{s_2y_k})\pi(t_{(s_3z_{i,l})^{-1}})=\sum\limits_{\substack{0\leq j\leq \min \{k, l\}}}\pi(t_{{}^\star(\hat z_{i,k+l-2j})})+\pi(t_{(\hat z_{i,k+l-2j})}).\end{aligned}$$

\medskip

The proof is completed.
\end{proof}

\medskip

{\bf Corollary 3.6I.}
For $x,y^{-1}\in\{s_2z_{i,j}, s_3z_{i,j}|\ i, j\ge 1\}$,  we have $\pi(t_xt_y)=\pi(t_x)\pi(t_y)$.
\begin{proof}
The result follows from Lemma 3.6G,  Lemma 3.6H, Lemma 3.6B, Lemma 3.6E, and Lemma 3.6F.
\end{proof}

\medskip

{ Combining Lemma 3.6F to Corollary 3.6I we have proved  assertion (b) in Proposition 3.6A. And assertion (c) follows by Corollary 3.5E.

Now we have proved Proposition 3.6A completely.}

\medskip

 {\bf Corollary 3.6K.}
 For any left cell $\Gamma\in Y_3,\ \Theta\in Y_1$,  we have a bijection
 $$\pi: \Gamma \cap\Theta^{-1}\to\{\text{iso. classes of irr. $F$-vector bundles on}\ \{\Theta,\ g\Theta\}\times \{\Gamma \}\}.$$
  The inverse gives a bijection
  $$\pi: \Theta\cap\Gamma ^{-1}\to\{\text{iso. classes of irr. $F$-vector bundles on}\ \{\Gamma \}\times \{\Theta,\ g\Theta\}\}.$$
\begin{proof}
It follows from Lemma 2.6 and subsection 2.4 and Proposition 3.6A.
\end{proof}

\bigskip

  {\bf 3.7.}  In this subsection we establish a bijection
  $$\pi: \Gamma_{02}\cap\Gamma_{013}^{-1}\to\{\text{iso. classes of irr. $F$-vector bundles on}\ Z_{013}\times \{\Gamma_{02}\}\}.$$
  The inverse gives a bijection
  $$\pi: \Gamma_{013}\cap\Gamma_{02}^{-1}\to\{\text{iso. classes of irr.$F$-vector bundles on}\ \{\Gamma_{02}\}\times Z_{013}\}.$$

By Lemma 3.5B, $\Gamma_{02}\cap\Gamma_{013}^{-1}\overset{\phi^{-1}}{\longrightarrow}
\Gamma_{02}\cap\Gamma_{03}^{-1}$ and $\Gamma_{013}^{-1}\overset{\phi^{-1}}{\longrightarrow}
\Gamma_{03}^{-1}$ are bijections,  and the set $\Gamma_{02}\cap\Gamma_{013}^{-1}=\{x\in W\ |\ \phi^{-1}(x)\in\Gamma_{02}\cap\Gamma_{03}^{-1}\}.$

\medskip

Denote the bijection in Proposition 3.6A as $\pi_4$.

{\bf Proposition 3.7A.} The map
$$\pi: \Gamma_{02}\cap \Gamma^{-1}_{013}\to\{\text{iso. classes of irr. $F$-vector bundles on}\ \{\Gamma_{02}\}\times Z_{013}
\}$$ defined by $\pi(x)=\pi_4\phi^{-1}(x)$ for any $x\in \Gamma_{02}\cap \Gamma^{-1}_{013}$ is a bijection.
\begin{proof}
It follows from Lemma 3.5D and Proposition 3.6A. \end{proof}

\medskip

  {\bf Corollary 3.7B.} For any left cell $\Gamma\in Y_3,\ \Theta\in Y_2$, we have a bijection
 $$\pi: \Gamma \cap\Theta^{-1}\to\{\text{iso. classes of irr. $F$-vector bundles on}\ \{\Theta,\ g\Theta\}\times \{\Gamma \}\}.$$
  The inverse gives a bijection
  $$\pi: \Theta\cap\Gamma ^{-1}\to\{\text{ iso. classes of irr. $F$-vector bundles  on}\ \{\Gamma \}\times \{\Theta,\ g\Theta\}\}.$$
\begin{proof}
It follows from Lemma 2.6 and subsection 2.4.
\end{proof}

\medskip

 {\bf 3.8.}  In this subsection we establish a bijection
  $$\pi: \Gamma_{03}\cap\Gamma_{013}^{-1}\to\{\text{iso. classes of irr. $F$-vector bundles on}\ Z_{013}\times Z_{03}\}.$$
  The inverse gives a bijection
  $$\pi: \Gamma_{013}\cap\Gamma_{03}^{-1}\to\{\text{iso. classes of irr. $F$-vector bundles on}\ Z_{03}\times Z_{013}\}.$$

By Lemma 3.5B, $\Gamma_{03}\cap\Gamma_{013}^{-1}\overset{\phi^{-1}}{\longrightarrow}
\Gamma_{03}\cap\Gamma_{03}^{-1}$ and $\Gamma_{013}^{-1}\overset{\phi^{-1}}{\longrightarrow}
\Gamma_{03}^{-1}$ are bijections,  and the set $\Gamma_{03}\cap\Gamma_{013}^{-1}=\{x\in W\ |\ \phi(x)\in\Gamma_{03}\cap\Gamma_{03}^{-1}\}.$

\medskip

Recall we denote the bijection in Proposition 3.4A as $\pi_2$.

{\bf Proposition 3.8A.}The map
$$\pi: \Gamma_{03}\cap \Gamma^{-1}_{013}\to\{\text{iso. classes of irr. $F$-vector bundles on}\ Z_{013}\times Z_{03}\}$$ defined by $\pi(x)=\pi_2\phi^{-1}(x)$ for any $x\in \Gamma_{03}\cap \Gamma^{-1}_{013}$ is a bijection.
\begin{proof}
It follows from Lemma 3.5D and Proposition 3.4A. \end{proof}

\medskip

  {\bf Corollary 3.8B.} For any left cell $\Gamma\in Y_1,\ \Theta\in Y_2$,  we have a bijection
 $$\pi: \Gamma \cap\Theta^{-1}\to\{\text{iso. classes of irr. $F$-vector bundles on}\ \{\Theta,\ g\Theta\}\times \{\Gamma,\ g\Gamma \}\}.$$
  The inverse gives a bijection
  $$\pi: \Theta\cap\Gamma ^{-1}\to\{\text{iso. classes of irr. $F$-vector bundles on}\ \{\Gamma, \ g\Gamma \}\times \{\Theta,\ g\Theta\}\}.$$
\begin{proof}
It follows from Lemma 2.6 and subsection 2.4.
\end{proof}

\bigskip

The proof for  Theorem 3.2 now is completed.

\bigskip

\noindent{\bf Acknowledgement:} Part of the work was done during
my visit to the Academy of Mathematics and Systems Science, Chinese Academy of Sciences. I am very grateful to the AMSS for hospitality and financial supports. And I would like to thank Nanhua Xi for useful discussions.



\begin{thebibliography}{99}

\bibitem[B]{B} R. Bezrukavnikov, {\sl On tensor categories attached to cells in affine Weyl groups,} In "Representation Theory of Algebraic Groups and Quantum Groups", Adv. Stud. Pure Math., 40, Math. Soc. Japan, Tokyo, 2004, pp. 69-90.



\bibitem[BO]{BO} R. Bezrukavnikov and V. Ostrik, {\sl On tensor categories attached to cells in affine Weyl groups II,} In "Representation Theory of Algebraic Groups and Quantum Groups", Adv. Stud. Pure Math., 40, Math. Soc. Japan, Tokyo, 2004,  pp.101-119.

\bibitem[DLP]{DLP} C. De Concini,  G. Lusztig,   C. Procesi,{\sl  Homology of the zero-set of a nilpotent vector field on a flag manifold}, J. Amer. Math. Soc. 1 (1988), 15-34.


\bibitem[D]{D} J. Du, {\sl The decomposition into cells of the affine Weyl group of type $\tilde{B_3}$,} Communications in Algebra, 16 (1988), no.7, 1383--1409.


\bibitem[KL]{KL} D. Kazhdan and G. Lusztig, {\sl Representations of Coxeter
groups and Hecke algebras,} Invent. Math. 53 (1979), 165-184.

\bibitem[L1]{L1} G. Lusztig, {\sl Cells in affine Weyl groups,} in
``Algebraic groups and related topics", Advanced Studies in Pure
Math., vol. 6, Kinokunia and North Holland, 1985, pp. 255-287.

\bibitem[L2]{L2} G. Lusztig, {\sl Cells in affine Weyl groups, II,} J. Alg.
109 (1987), 536-548.

\bibitem[L3]{L3} G. Lusztig, {\sl Cells in affine Weyl groups, IV,} Journal of The Faculty of Science, 36 (1989), no.2, 297-328.



\bibitem[X1]{X1} N. Xi, {\sl Representations of Affine Hecke Algebras,} volume 1587, Springer Lecture Notes in Math., 1994.

\bibitem[X2]{X2} N. Xi, {\sl The based ring of two-sided cells of affine Weyl groups of type ${\tilde{A}_{n-1}}$,} volume 749, American Mathematical Soc., 2002.

 \bibitem[QX]{QX} Y.Qiu and N. Xi, {\sl The based ring of two-sided cells in an affine Weyl group of type $\tilde B_3$, I,} Sci. China Math., 2022.

 \bibitem[QX2]{QX2} Y. Qiu and N. Xi, {\sl The based ring of two-sided cells in an affine Weyl group of type $\tilde B_3$, II,} Pure and Applied Mathematics Quarterly, to appear, arxiv: 2202.00302.



\end{thebibliography}
\end{document}